\documentclass[ejs]{imsart2} % for arxiv version

\RequirePackage[OT1]{fontenc}
\RequirePackage{amsthm,amsmath}
\RequirePackage[numbers,sort&compress]{natbib}
\RequirePackage[colorlinks,citecolor=blue,urlcolor=blue]{hyperref}
\RequirePackage{amsfonts,amsmath}
\RequirePackage{graphicx}
\RequirePackage{algorithmic}
\RequirePackage{algorithm}
\RequirePackage{comment}
\RequirePackage[document]{ragged2e}

\pubyear{2017}
\volume{0}
\issue{0}
\firstpage{0}
\lastpage{0}

\startlocaldefs
\numberwithin{equation}{section}
\theoremstyle{plain}
\newtheorem{theorem}{Theorem}[section]
\endlocaldefs

\newcommand{\ber}{\begin{eqnarray}}
\newcommand{\eer}{\end{eqnarray}}

\newtheorem{lemma}{\noindent Lemma}
\newtheorem{corollary}{\noindent Corollary}

\newtheorem{remark}{\noindent Remark}

\newcommand{\be}{\begin{equation}}
\newcommand{\ee}{\end{equation}}
\newcommand{\bal}{\begin{align}}
\newcommand{\eal}{\end{align}}
\newcommand{\balnonum}{\begin{align*}}
\newcommand{\ealnonum}{\end{align*}}

\newcommand{\hygeom}{{}_2 F_1}

\newcommand{\RS}{Robbins and Siegmund\:}
\newcommand{\DD}{Dudewicz and Dalal\:}

\def\qed{\hfill \vrule height1.3ex width1.2ex depth-0.1ex}

\begin{document}

\begin{frontmatter}
\title{On the Asymptotic Efficiency of Selection Procedures for Independent Gaussian Populations}
\runtitle{Asymptotic Efficiency of Selection Procedures.}

\begin{aug}
\author{\fnms{Royi} \snm{Jacobovic}\ead[label=e1]{royi.jacobovic@mail.huji.ac.il}}
\and
\author{\fnms{Or} \snm{Zuk}\ead[label=e2]{or.zuk@mail.huji.ac.il }}

\address{Department of Statistics \\ The Hebrew University of Jerusalem \\ Mt.Scopus, Jerusalem, 91905.\\
\printead{e1,e2}}

\runauthor{R. Jacobovic and O. Zuk}

\end{aug}

\begin{abstract}
The field of discrete event simulation and optimization techniques motivates researchers to adjust classic ranking and selection (R\&S) procedures to the settings where the number of populations is large. We use insights from extreme value theory in order to reveal the asymptotic properties of R\&S procedures. Namely, we generalize the asymptotic result of \RS regarding selection from independent Gaussian populations with known constant variance by their means to the case of selecting a subset of varying size out of a given set of populations. In addition, we revisit the problem of selecting the population with the highest mean among independent Gaussian populations with unknown and possibly different variances. Particularly, we derive the relative asymptotic efficiency of \DD's and Rinott's procedures, showing that the former can be asymptotically superior by a multiplicative factor which is larger than one, but this factor may be reduced by proper choice of parameters. We also use our asymptotic results to suggest that the sample size in the first stage of the two procedures should be logarithmic in the number of populations. 
\end{abstract}



\end{frontmatter}

\section{Introduction}

Selecting and ranking items from a set based on incomplete and noisy information is a natural problem arising in many domains with limited resources. Examples include selecting students for a program from a list of candidates based on their prior grades, ranking web-pages based on their relevance to a query and displaying the top pages to a user, or 
finding the best (or near the best) system design with respect to some measure of performance. Discrete event simulation is a popular methodology for studying such system design problems, with some reviews of applications in \cite{forgionne1983corporate,harpell1989operations,lane1993operations,shannon1980operation}. Fundamental texts summarizing the basics of this approach are \cite{fu1994optimization} and \cite{schruben1989review}. More general references for stochastic simulations are given in \cite{banks1998handbook,chen2011stochastic,fu2015handbook,fu2005simulation,kleijnen2008design}.

The modern literature about discrete event simulation is strongly related to the theory of ranking and selection (R\&S) procedures. This literature considers a set of populations and a user who wants to select the populations associated with a specific relative stochastic property such as the highest mean, the smallest variance, etc. With regard to this task, the R\&S literature is devoted to development of useful procedures, i.e. sampling policies and selection rules to pinpoint the target populations with some performance guarantee and low sampling cost. A nice glance into the R\&S theory is provided in \cite{bechhofer1990comparison} while extensive summary can be found in the books \cite{gibbons1999selecting,gupta2002multiple}. Major fields of this research include Bayesian and indifference-zone (IZ) formulations. Recent work in the Bayesian context is summarized by \cite{Chick2006,chick2006bayesian}, and an extension to the case of multiple attributes appears in \cite{Frazier2011}. Similarly, recent contributions regarding the IZ formulation are described in \cite{kim2006a}.

As demonstrated by \cite{goldsman1998comparing}, the R\&S literature offers attractive procedures for the case where the number of alternative designs is relatively small and there is no strong functional relationship among them. However, as pointed by \cite{kim2006}, this situation is not frequent in practice. In particular, the number of alternative designs is usually large which means that classical R\&S procedures cannot be applied directly with no proper adjustments. Motivated by this issue, several authors introduced improvements of the classic R\&S procedures as a solution for this problem \cite{ahmed2002simulation,boesel2003using,chen2000,nelson2001simple,Frazier2014}. These improvements were mostly compared to their classic R\&S ancestors by simulations (although rigorous bounds were derived in \cite{Frazier2014} for the fully sequential case).
While simulations can be carried out to study these modern procedures, they do not provide insights or rigorous bounds regarding the performance as a function of procedures' choices and parameters, and become computationally intensive as the number of populations and sample size grow. A complementary and attractive approach is to evaluate the quality of R\&S procedures by investigating their asymptotic behavior, with the goal being providing rigorous analytic bounds and approximations for the procedures' performance, thus gaining insights into their dependence on various parameters and on the relative efficiency of different procedures. The fundamentals of the asymptotic theory of R\&S procedures appear in the book \cite{mukhopadhyay1994multistage}. This work makes more contributions to this theory.

In details, Section \ref{sec:linear_maxima} applies insights from extreme value theory to specify the asymptotic behavior of linear combinations of maxima. Sections \ref{sec:robbins_siegmund} and \ref{sec:two_stage} use these results in order to derive new asymptotic results for well-known R\&S procedures through the IZ approach of Bechhofer \cite{bechhofer1954single}. Namely, Section \ref{sec:robbins_siegmund} generalizes the result of \RS \cite{robbins1967iterated} who considered the problem of selection from $k$ independent normal homoscedastic populations with known variance by their means. \RS provided a first order approximation for the minimal sample-size which controls the probability for correct selection ({\it PCS}) of the {\it single} population with highest mean as the total number of populations tends to infinity \cite{robbins1967iterated}. This work generalizes their results to the case where the number of selected populations can be determined as a function of the total number of populations, and deriving the asymptotic sample size required to achieve a desired {\it PCS} as a function of the number of selected populations. In addition, we present a new proof for the original result. Section \ref{sec:two_stage} starts by brief review of two well-known two-stage procedures which were proposed respectively by \DD \cite{dudewicz1975allocation} and Rinott \cite{rinott1978two}. Both procedures were designed for the problem of selecting the Gaussian population with the highest mean for independent populations with unknown and possibly different variances. We derive first order approximations for these procedures asymptotic efficiencies, measured in terms of the expected sample size required to achieve a desired {\it PCS}, as the total number of populations grows to infinity. A corollary of these results is that asymptotically, Rinott's procedure is relatively less efficient than the procedure of \DD by a multiplicative factor depending on the initial sample size used in stage one of both procedures. However, our asymptotic analysis motivates a conjecture that the optimal sample size in the first stage of both procedures grows logarithmically in the number of populations, and with this optimal choice the multiplicative factor approaches one and the two procedures may be asymptotically equivalent.

We performed numerical computations in order to highlight and complement our analytic asymptotic results - {\it Matlab} code for these computations, including a script reproducing all figures in the paper is available from {\it github} at \\ {\small \url{https://github.com/orzuk/MatUtils/tree/master/stats/ranking_selection} .} 

\newpage

\section{Asymptotics of Linear Combinations of Partial Maxima}
\label{sec:linear_maxima}
Let $\{X_m;m=1,2,\ldots\}$ be an infinite sequence of identically independently distributed (i.i.d) continuous random-variables (r.v's) with cumulative distribution function (c.d.f) $F$ such that $F(-\infty)=0$ and $\underset{x\rightarrow\infty}{\lim} F(x) =1$. Let $T\in\mathbb{N}$ and for each $k>T$ let $1=:s_k^{(0)} < s_k^{(1)}<\ldots<s_k^{(T-1)}<s_k^{(T)}:=k$ be an increasing integer sequence defining a partition of $\{1,..,k\}$ into $T$ sub-groups.
Define the partial maxima of $\{X_m;m=1,\ldots,k\}$ with respect to this partition by $M_k^{(t)}:=\max\{X_j;s_k^{(t-1)}+1\leq j\leq s_k^{(t)}\}$, $t\in\mathcal{T}:=\{1,\ldots,T\}$. With regard to this sequence of partitions assume that for each $t\in\mathcal{T}$ the difference $\delta_k^{(t)} :=s_k^{(t)}-s_k^{(t-1)}$ converges in the broad sense, i.e. there exist $\delta^{(t)}\in\mathbb{N}\cup\{\infty\}$ such that $\delta_k^{(t)} \overset{k \to \infty}{\longrightarrow} \delta^{(t)}$. Moreover, let $F$ be max-stable in the sense that it is associated with an extreme value distribution, i.e. there are two sequences of normalizing constants $\{a_k\}$ and $\{b_k\}$ such that:

\begin{enumerate}
\item There exists $K\in\mathbb{N}$ such that $a_k>0$ for any $k>K$.
\item $\{a_k\}$ is weakly-monotonic.
\item $a_k\big(\underset{j=1,\ldots,k}{\max} X_j-b_k\big) \xrightarrow {\mathcal{L}}Y$ where the notation $\xrightarrow{\mathcal{L}}$ denotes convergence in law of r.v's and $Y\sim F_Y$ is a continuous r.v, i.e. its c.d.f $F_Y$ is characterized by Fisher-Tippet-Gnedenko's theorem.
\end{enumerate}

Considering the deterministic sequence $\{\xi_k;k\in\mathbb{N}\}\subset\mathbb{R}$ and some vector $\alpha:=(\alpha_1,\ldots,\alpha_T)\in\mathbb{R}^T$, the goal of this section is to calculate the following limit:

\be
L:=L\Big(\{\delta_k^{(1)}\},\ldots,\{\delta_k^{(T)}\},\{\xi_k\},\alpha; F, F_Y\Big)=\lim\limits_{k\rightarrow\infty}\mathbb{P}\big(\sum_{t=1}^{T}\alpha_t M_k^{(t)}\leq\xi_k\big) .
\ee	

To phrase the main results, consider the partition $\mathcal{T}=\mathcal{T}_1\cup\mathcal{T}_2$ defined by the sets $\mathcal{T}_2:=\{t\in\mathcal{T};\delta^{(t)}=\infty\}$, corresponding to infinite subsequences, and $\mathcal{T}_1:=\mathcal{T}\setminus\mathcal{T}_2$, corresponding to finite subsequences. For $t\in\mathcal{T}_1$, the limit $\lim\limits_{k\rightarrow\infty} a_{\delta_k^{(t)}} = a_{\delta^{(t)}}\in\mathbb{R}_{++}$ exists. In addition, since $\{a_k\}$ is positive and weakly-monotonic, the limit $a_\infty:=\lim\limits_{k\rightarrow\infty}a_k\in[0,\infty]$ exists in the broad sense. With regard to this framework, our main theorems provide sufficient conditions under which $L$ exists and can be calculated:

\begin{theorem}
\label{thm:sum_max_a_k_positive}
Let $F$ be a max-stable distribution associated with sequences 
$\{a_k\}, \{b_k\}$ such that $a_\infty\in(0,\infty]$. For each $t\in\mathcal{T}_1$,
let $M_t$ be a random variable distributed as the maximum of $\delta^{(t)}$ i.i.d. random variables with c.d.f. $F$, i.e. $M_t \sim F^{\delta^{(t)}}$ and for each $t\in\mathcal{T}_2$, let $Y_t\sim F_Y$ such that $\{M_t;t\in\mathcal{T}_1\}\cup\{Y_t;t\in\mathcal{T}_2\}$ is a set of independent r.v's. Define $V$ as		
\be
V := \sum_{t\in\mathcal{T}_1}\alpha_tM_t+\sum_{t\in\mathcal{T}_2}\frac{\alpha_t}{a_\infty}Y_t
\ee				
and suppose that the limit 		
\be
L^*:=\lim\limits_{k\rightarrow\infty}\bigg\{\xi_k-\sum_{t\in\mathcal{T}_2}\alpha_t b_{\delta_k^{(t)}}\bigg\}
\ee
exists in the broad sense, i.e. $L^*\in\bar{\mathbb{R}}$. Then: 
\be
L=F_V(L^*).
\ee
\end{theorem}

\begin{theorem}
\label{thm:sum_max_a_k_zero}
Suppose that $F$ is associated with $\{a_k\},\{b_k\}$ such that $a_{\infty}=0$, 
there exists an index $t^*\in\mathcal{T}_2$ such that			
\be
\lambda_t:=\lim\limits_{k\rightarrow\infty}\frac{a_{\delta_k^{(t^*)}}}{a_{\delta_k^{(t)}}}\in\mathbb{R}_+\ , \ \forall t\in\mathcal{T}_2 \:, 
\ee
$V$ is given by		
\be
V:=\sum_{t\in\mathcal{T}_2} \alpha_t \lambda_t Y_t \ \ s.t. \ \ Y_t\stackrel{i.i.d.}{\sim} F_Y,t\in\mathcal{T}_2
\ee		
and the limit 		
\be		
L^{**}:=\lim\limits_{k\rightarrow\infty}a_k(\xi_k-\sum_{t\in\mathcal{T}_2}\alpha_t b_{\delta_k^{(t)}})
\ee		
exists in the broad sense, i.e. $L^{**}\in\bar{\mathbb{R}}$. Then: 
\be
L=F_V(L^{**}).
\ee		
\end{theorem}

\subsection{Proofs}

The proofs of Theorems \ref{thm:sum_max_a_k_positive} and \ref{thm:sum_max_a_k_zero} are based on the following two lemmata on convergence in law:	
\begin{lemma}
\label{lemma:convergence_law}
If $V$ is a finite r.v and $V_1,V_2\ldots$ are r.v's such that

\begin{enumerate}
\item $V_k\xrightarrow{\mathcal{L}}V$.

\item $F_V$ is continuous on $\mathcal{C}\subseteq\mathbb{R}$.

\item $\{x_k;k\in\mathbb{N}\}$ is a deterministic sequence such that $x_k\rightarrow\bar{x}$ where $\bar{x}\in\mathcal{C}\cup\{-\infty,\infty\}$ .
\end{enumerate} 

Then $F_{V_k}(x_k)\xrightarrow[k\rightarrow\infty]{} F_V(\bar{x})$.
\end{lemma}

\begin{lemma}
\label{lemma:convergence_sum}
Let $s \in \mathbb{N}$. If $\forall k\in\mathbb{N}$, $Z_k^{(1)},\ldots,Z_k^{(s)}$ are independent r.v's such that $Z_k^{(i)}\xrightarrow[k\rightarrow\infty]{\mathcal{L}}Z_i, \: \forall i=1,\ldots,s$, then
\be		
Z_k^{(1)}+\ldots +Z_k^{(s)}\xrightarrow[k\rightarrow\infty]{\mathcal{L}}Z_1+\ldots+Z_s
\ee		
where $Z_1,\ldots,Z_s$ are independent r.v's.
\end{lemma}
Lemma \ref{lemma:convergence_law} is a known result about convergence in law. More details are provided in \cite{billingsley2013convergence}. Lemma \ref{lemma:convergence_sum} is obtained by a straightforward application of the multivariate continuous mapping theorem for the vector $(Z_k^{(1)},\ldots,Z_k^{(s)})$, noticing that due to independence we have 
$(Z_k^{(1)},\ldots,Z_k^{(s)})\xrightarrow{\mathcal{L}}(Z_1,\ldots,Z_s)$.

\proof(Theorem \ref{thm:sum_max_a_k_positive})

Assume first that $\frac{\alpha_t}{a_{\infty}}\neq0, \: \forall t\in\mathcal{T}$ and express the limit $L$ as follows: 	
\begin{align}	L &= \lim\limits_{k\rightarrow\infty}\mathbb{P}\big(\sum_{t=1}^{T}\alpha_t M_k^{(t)}\leq\xi_k\big) \nonumber \\
&=\lim\limits_{k\rightarrow\infty}\mathbb{P}\big(\sum_{t\in\mathcal{T}_1}\alpha_t M_k^{(t)}+\sum_{t\in\mathcal{T}_2}\alpha_t M_k^{(t)}\leq\xi_k\big) \nonumber \\
&= \lim\limits_{k\rightarrow\infty}\mathbb{P}\big(\sum_{t\in\mathcal{T}_1}\alpha_t M_k^{(t)}+\sum_{t\in\mathcal{T}_2}\frac{\alpha_t}{a_{\delta_k^{(t)}}} a_{\delta_k^{(t)}}(M_k^{(t)}-b_{\delta_k^{(t)}})\leq\xi_k-\sum_{t\in\mathcal{T}_2}\alpha_t b_{\delta_k^{(t)}}\big) .
\end{align}	
For any $t\in\mathcal{T}_2$, known properties of convergence in law imply that 	
\be	
\frac{\alpha_t}{a_{\delta_k^{(t)}}} a_{\delta_k^{(t)}}(M_k^{(t)}-b_{\delta_k^{(t)}})\xrightarrow{\mathcal{L}}\frac{\alpha_t}{a_\infty}Y_t
\ee	
In addition, $\forall t\in\mathcal{T}_1$, $\alpha_t M_k^{(t)}\xrightarrow{pointwise}\alpha_t M_t$ and hence $\alpha_t M_k^{(t)}\xrightarrow{\mathcal{L}}\alpha_t M_t$. For any $k\in\mathbb{N}$, the random variables $M_k^{(1)},\ldots,M_k^{(T)}$ are determined by disjoint subgroups of i.i.d sequence of r.v's and consequently $\forall k\in\mathbb{N}$, $M_k^{(1)},\ldots,M_k^{(T)}$ are independent r.v's. Therefore, Lemma \ref{lemma:convergence_sum} implies that	
\be
\sum_{t\in\mathcal{T}_1}\alpha_t M_k^{(t)}+\sum_{t\in\mathcal{T}_2}\frac{\alpha_t}{a_{\delta_k^{(t)}}} a_{\delta_k^{(t)}}(M_k^{(t)}-b_{\delta_k^{(t)}})\xrightarrow{\mathcal{L}}\sum_{t\in\mathcal{T}_1}\alpha_t M_t+\sum_{t\in\mathcal{T}_2}\frac{\alpha_t}{a_\infty} Y_t=V
\label{eq:sum_two_parts}
\ee	
where $\{M_t;t\in\mathcal{T}_1\}\cup\{Y_t;t\in\mathcal{T}_2\}$ is a set of independent r.v's. 

At this stage, assume that $\exists t'\in\mathcal{T}_2$ for which $\frac{\alpha_{t'}}{a_\infty}=0$. The LHS of eq. (\ref{eq:sum_two_parts}) can be represented as the sum of two finite sums $\mathcal{S}_1+\mathcal{S}_2$, where $\mathcal{S}_1$ includes all summands that converge in law to zero and $\mathcal{S}_2$ includes all other summands. Recalling that convergence in law to a constant implies convergence in probability, then each of the summands in $\mathcal{S}_1$ converges in probability to zero, and since the number of summands is finite, $\mathcal{S}_1 \overset{\mathbb{P}}{\rightarrow} 0$ . Similarly, by the arguments used under the simplifying assumption that $\frac{\alpha_t}{a_\infty}\neq0,\forall t\in\mathcal{T}_2$, $\mathcal{S}_2 \overset{\mathbb{L}}{\rightarrow} V$. Therefore, by Slutsky's Lemma (see Chapter $6$ in \cite{ferguson1996course}), the total sum converges in law to $V$. The distribution of $Y$ is characterized by Fisher-Tippet-Gnedenko's Theorem hence $Y$ is a continuous r.v. In addition, $M_t$ is distributed like a maximum of a finite number of i.i.d continuous r.v's, and is a continuous r.v. Therefore, deduce that $V$ is a finite sum of independent continuous r.v's and hence it is a continuous r.v. Finally, since $L^*$ exists in the broad sense, the needed result follows directly from Lemma \ref{lemma:convergence_law}. 	
\endproof

\proof(Theorem \ref{thm:sum_max_a_k_zero})
In the spirit of the proof of Theorem \ref{thm:sum_max_a_k_positive}, it is enough to prove the theorem under the simplifying assumption $\alpha_t\lambda_t\neq 0, \: \forall t\in\mathcal{T}_2$. Under this assumption, the limit $L$ can be expressed as follows: 	
\begin{align}	L &= \lim\limits_{k\rightarrow\infty}\mathbb{P}\big(\sum_{t=1}^{T}\alpha_t M_k^{(t)}\leq\xi_k\big) \nonumber \\
&= \mathbb{P}\big(\sum_{t\in\mathcal{T}_1}\alpha_t M_k^{(t)}+\sum_{t\in\mathcal{T}_2}\alpha_t M_k^{(t)}\leq\xi_k\big) \nonumber \\
&= \lim\limits_{k\rightarrow\infty}\mathbb{P}\Big(\sum_{t\in\mathcal{T}_1}\alpha_t a_{\delta^{t^*}_k}M_k^{(t)}+\sum_{t\in\mathcal{T}_2}\alpha_t\frac{a_{\delta^{t^*}_k}}{a_{\delta_k^{(t)}}} a_{\delta_k^{(t)}} (M_k^{(t)}-b_{\delta_k^{(t)}})\leq a_{\delta^{t^*}_k}\big(\xi_k-\sum_{t\in\mathcal{T}_2}\alpha_tb_{\delta_k^{(t)}}\big)\Big) . 
\end{align}	
By similar arguments as in the proof of Theorem \ref{thm:sum_max_a_k_positive}, 	
\begin{enumerate}
\item $\alpha_t a_{\delta^{t^*}_k}M_k^{(t)}\xrightarrow{\mathcal{L}}0$ , $\forall t\in\mathcal{T}_1$.
\item $\alpha_t \frac{a_{\delta^{t^*}_k}}{a_{\delta_k^{(t)}}} a_{\delta_k^{(t)}} (M_k^{(t)}-b_{\delta_k^{(t)}})\xrightarrow{\mathcal{L}}\alpha_t \lambda_tY_t$, $\forall t\in\mathcal{T}_2$.
\item For any $k\in\mathbb{N}$, $M_k^{(1)},\ldots,M_k^{(T)}$ are independent r.v's.
\end{enumerate}

Since all the preconditions of Lemma \ref{lemma:convergence_sum} are satisfied, deduce that 	
\be		
\sum_{t\in\mathcal{T}_2}\alpha_t \frac{a_{\delta^{t^*}_k}}{a_{\delta_k^{(t)}}} a_{\delta_k^{(t)}} (M_k^{(t)}-b_{\delta_k^{(t)}})\xrightarrow{\mathcal{L}}\sum_{t\in\mathcal{T}_2}\alpha_t \lambda_t Y_t=V\ \ .
\ee		
where $\{Y_t;t\in\mathcal{T}_1\}$ are independent r.v's. On the other hand, because convergence in law to a constant implies convergence in probability,
\be			
\sum_{t\in\mathcal{T}_1}\alpha_t a_{\delta^{t^*}_k}M_k^{(t)}\overset{\mathbb{P}}{\longrightarrow} 0 . 
\ee		
Thus, Slutsky's Lemma can be applied to obtain the following limit:
\be		\sum_{t\in\mathcal{T}_1}\alpha_t a_{\delta^{t^*}_k}M_k^{(t)}+\sum_{t\in\mathcal{T}_2}\alpha_t\frac{a_{\delta^{t^*}_k}}{a_{\delta_k^{(t)}}} a_{\delta_k^{(t)}} (M_k^{(t)}-b_{\delta_k^{(t)}})\overset{\mathcal{L}}{\longrightarrow} 0+V=V .
\ee		
Fisher-Tippet-Gnedenko's Theorem implies that $\{Y_t;t\in\mathcal{T}_1\}$ are continuous r.v's. Therefore, $V$ is a finite sum of finite continuous independent r.v's, so it is also a finite continuous r.v. Finally, since $L^*$ exists in the broad sense, Lemma \ref{lemma:convergence_law} implies the needed result.		
\endproof

\section{Generalized Robbins-Siegmund Result}
\label{sec:robbins_siegmund}
This Section demonstrates an application of Theorem \ref{thm:sum_max_a_k_positive} to the problem of selection from homoscedastic independent Gaussian populations with known variance by their means. Subsection \ref{sec:stat_framework} depicts the relevant statistical model. Subsection \ref{sec:asymptotic_sample_size}
includes a short review of the original result \cite{robbins1967iterated} as well as our generalization of this result to the case of selecting more than one population.

\subsection{Statistical Framework}
\label{sec:stat_framework}
Let $X_{ij}\sim N(\theta_i,\sigma^2);i=1,\ldots,k, \: j=1,\ldots,N$ be independent univariate Gaussian r.v's with known variance $\sigma^2>0$ and unknown means $\theta=(\theta_1,\ldots,\theta_k)\in\mathbb{R}^k$. The task is to find the $1\leq s\leq \lfloor \frac{k}{2}\rfloor$ populations with the largest means. An intuitive procedure for this purpose is to compute the empirical means $\bar{X}_i:=\frac{1}{N}\sum_{j=1}^NX_{ij}, \: \forall i=1,\ldots,k$ and select the $s$ populations associated with the highest values. This procedure can be justified theoretically as explained in Chapter $3$ of \cite{gupta2002multiple}. Our goal is to find the minimal $N$ which ensures correct selection of all of the required $s$ populations with probability of at least $p\in(0,1)$ for this procedure.

It is not possible to control the probability of correct selection (denoted by $\mathbb{P}(CS)$) without further assumptions. To see this, take populations with equal means, $\theta_{\gamma} := \gamma \cdot \textbf{1}_k;\gamma\in\mathbb{R}$. There is no way to distinguish between the populations by sampling from them and the $\mathbb{P}(CS)$ is not sensitive to $N$. In order to allow the user to distinguish between populations, the indifference-zone approach of Bechhofer \cite{bechhofer1954single} is adopted, i.e. the parameter-space is restricted in the following way	
\be
\theta \in \Theta(\Delta,k) := \{\tilde{\theta}\in\mathbb{R}^k;\tilde{\theta}_{[k-s+1]}-\tilde{\theta}_{[k-s]}\geq\Delta\}
\label{eq:indifference_parameter_space}
\ee	
where $\tilde{\theta}_{[1]}\leq\ldots\leq\tilde{\theta}_{[k]}$ are the ordered means and $\Delta>0$ is a known parameter indicating the minimal difference in mean between the top $s$ and bottom $k-s$ populations. $\Delta$ can also be interpreted as the indifference level of the experimenter, i.e. if the absolute value of the difference between the means of two different population is less than $\Delta$, the experimenter will consider them as equivalent populations. Since $(\bar{X}_1,\ldots,\bar{X}_k)$ is a consistent estimator for any $\theta\in\Theta(\Delta,k)$, the probability of correct selection tends to $1$ as $N\rightarrow\infty$, regardless of the true parametrization. Thus, the minimal sample size which ensures correct selection with probability $p$ is well-defined. 

\subsection{Asymptotic Sample-Size}
\label{sec:asymptotic_sample_size}	
Fix $\Delta>0$ and denote by $N^*_{k,s}(p)$ the minimal $N\in\mathbb{N}$ for which the probability of correct selection is above $p$. For simplicity, we follow \cite{gupta2002multiple} and ignore the rounding error, i.e. $N_{k,s}^*(p)$ is defined as the solution of the following equation
\be
\mathbb{P}(CS_{k,N_{k,s}^*}^s;\Delta,\theta^*)=p
\label{eq:h_1_integral}
\ee	
where $CS_{k,N}^s$ is the event of making correct selection of the $s$ out of $k$ populations with the highest means based on $N$ samples from each population. In addition, $\theta^*$ is some \textit{least favorable configuration (LFC)}, i.e. it is a parametrization which satisfies 	
\be
\mathbb{P}(CS^s_{k,N_{k,s}};\Delta, \theta^*)=\inf_{\tilde{\theta}\in\Theta(\Delta,k)}\mathbb{P}(CS^s_{k,N_{k,s}};\Delta,\tilde{\theta})\ \ .
\ee	

With regard to this model, it was shown in \cite{robbins1967iterated} that for $s=1$ and $k\rightarrow\infty$, $N_{k,s=1}^*(p)\sim\frac{2\sigma^2}{\Delta^2}\ln(k-1)$ regardless of the value of $p\in(0,1)$. The asymptotic notation $\sim$ is interpreted in its classical terminology, i.e. for any two sequences $a_k\sim b_k$ if and only if $\frac{a_k}{b_k} \overset{k \to \infty}{\longrightarrow} 1$. 

We study the more general settings where the inequality $1\leq s:=s_k<k-s_k$ holds asymptotically and both limits $s_k\rightarrow\bar{s}\in\mathbb{N}\cup\{\infty\}$ and $\frac{\ln(s_k)}{\ln(k-s_k)}\rightarrow C\in[0,1]$ exist (e.g. $s_k\equiv \bar{s}\in\mathbb{N}$, $s_k:=\ln(k)$ and $s_k:=\psi k^{\beta}; 0<\psi,0<\beta \leq 1$ with $\psi<\frac{1}{2}$ for $\beta=1$). To simplify notation, let $N^*_k(p):=N^*_{s_k,k}(p)$ when $s_k$ is understood from context. Theorem \ref{thm:generalized_sigmound_robbins} shows that for each combination of $p\in(0,1)$ and $k \gg 1$, $N^*_k(p)$ is well defined and gives the asymptotics of $N^*_k(p)$ as $k \to \infty$, showing that the first order does not depend on $p$. Observe that in the following statement we can replace $\ln(k-s_k)$ by $\ln(k)$ in the asymptotic expression $\tilde{N}_k^*$ because the assumption that asymptotically $1\leq s_k\leq\frac{k}{2}$ implies $\ln(k-s_k)\sim\ln(k)$.

\begin{theorem}
\label{thm:generalized_sigmound_robbins}		
Let $s_k$ be a sequence such that: (i) $\forall k\in\mathbb{N}$, $1\leq s_k\leq k-s_k$.
(ii) $s_k\rightarrow\bar{s}\in\mathbb{N}\cup\{\infty\}$, and 
(iii) $\frac{\ln(s_k)}{\ln(k-s_k)}\rightarrow C\in[0,1]$. 
Then, $\forall p\in(0,1)$ the following statements hold:

\begin{enumerate}
\item There exists $K_p\in\mathbb{N}$ such that $N^*_k(p)$ exists for any $k>K_p$.
\item $N^*_k(p) \sim \tilde{N}_k^* := \frac{2\sigma^2(1+\sqrt{C})^2}{\Delta^2}\ln(k-s_k)$ as $k\rightarrow\infty$.			
\end{enumerate} 	
\end{theorem}

\proof

\begin{enumerate}
\item
Recall the definition of $\Theta(\Delta,k)$ and the fact that the $X_{ij}$'s are Gaussians with equal variances. Under these assumptions, the collection of all of the least favorable configurations is given by:
\be
\Theta^*(\Delta,k):=\{\theta^c\in\mathbb{R}^k;\: \theta^c_{[i]}=c+\Delta \cdot \textbf{1}_{\{i>k-s_k\}};\: c\in\mathbb{R}\} . 
\label{eq:lfc}
\ee	

where $\textbf{1}_A$ denotes the standard indicator function returning $1$ for any input $x \in A$ and $0$ otherwise. Consequently, w.l.o.g. set $c=0$, $Z=(Z_1,\ldots,Z_k)\sim N_k(0,I)$ and define the function 	
\begin{align}
f(n,k) &:= \inf_{\tilde{\theta}\in\Theta(\Delta,k)}\mathbb{P}(CS^{s_k}_{k,n};\Delta,\tilde{\theta})=P(CS^{s_k}_{k,n};\theta^0) \nonumber \\
&= \mathbb{P}(\bar{X}_i<\bar{X}_j,\forall 1\leq i \leq k-s_k < j\leq k) \nonumber \\
&= \mathbb{P}(\max_{i=1,\ldots ,k-s_k}\bar{X_i}\leq\min_{j=k-s_k+1,\ldots ,k}\bar{X_j}) \nonumber \\
&= \mathbb{P}\Big(\frac{\sigma}{\sqrt{n}}\max_{i=1,\ldots ,k-s_k+1}Z_i\leq\frac{\sigma}{\sqrt{n}}\min_{j=k-s_k+1,\ldots ,k}Z_j+\Delta\Big) . 
\end{align}
Denote the partial maxima by $M_k^{(1)} := \underset{i=1,\ldots,k-s_k}{\max}Z_i$ ; $M_k^{(2)} := \underset{i=k-s_k+1,\ldots,k}{\max}Z_i$. 
Since the univariate centered Gaussian is symmetric around zero, 
it is possible to express $f(n,k)$ in terms of a linear combination of $M_k^{(1)}$ and $M_k^{(2)}$. Since $M_k^{(1)},M_k^{(2)}$ are random variables determined by disjoint subsets of $Z_1,\ldots,Z_k$, they are independent. Therefore, $f(n,k)$ can be expressed by the following convolution in terms of the p.d.f $\phi$ and c.d.f $\Phi$ of the standard Gaussian distribution,	
\be
f(n,k)=\mathbb{P}\Big[M_k^{(2)}+M_k^{(1)}\leq\frac{\Delta\sqrt{n}}{\sigma}\Big]=
\int_{-\infty}^{\infty}\Phi^{k-s_k}\big(\frac{\Delta\sqrt{n}}{\sigma}-t\big)s_k\Phi^{s_k-1}(t)\phi(t)dt .
\ee	
Using this representation of $f(n,k)$ as a continuous c.d.f. for fixed $k$, we observe that:

\begin{enumerate}
\item $f(\infty;k) := \underset{n \to \infty}{\lim} f(n;k) = 1, \: \forall k\in\mathbb{N}$.

\item For fixed $k$ and $\forall n>0$, $f_k(n):=f(n;k)$ is continuous in $n$. 

\item As mentioned in \cite{leadbetter2012extremes}, page $20$, Example $1.7.1$, the normalizing sequences of the maximum of standard Gaussian r.v's are given by:
\begin{align}
a_k &= \sqrt{2\ln(k)}\ , \ \ \forall k\in\mathbb{N} \nonumber \\
b_k &= \sqrt{2\ln(k)}-\frac{\ln\ln(k)-\ln(4\pi)}{2\sqrt{2\ln(k)}}\ , \ \ \forall k\in\mathbb{N} .
\end{align}		
Since the univariate Gaussian is a continuous r.v and $a_k\rightarrow\infty$, $\underset{k \to \infty}{\lim} f(1,k)$ can be phrased in the form described by Theorem \ref{thm:sum_max_a_k_positive}: set $\alpha=(1,1)$, $\xi_k\equiv\frac{\Delta}{\sigma}$ and define $V$ according to the description of Theorem \ref{thm:sum_max_a_k_positive}. Observe the following limit		
\be		
L^*=\lim\limits_{k\rightarrow\infty}\big\{\frac{\Delta}{\sigma}-b_{k-s_k}-b_{s_k}\big\}=-\infty , 
\ee		
i.e. Theorem \ref{thm:sum_max_a_k_positive} implies $\underset{k \to \infty}{\lim} f(1,k) = 0$. 
\end{enumerate}

To end the proof of the first statement, fix $p\in(0,1)$. By Observation $3$, $\exists K_p\in\mathbb{N}$ such that $f_{k}(1)=f(1,k)<p, \: \forall k>K_p$. Let $k>K_p$. According to Observation $1$, there is $\bar{n}>1$ such that $f_{k}(\bar{n})>p$. Since Observation $2$ states that $f_{k}$ is continuous on $[1,\bar{n}] \subset (0,\infty)$, by the intermediate-value theorem $\exists n^*>1$ with $f(n^*,k)=p$. Thus, $\forall k>K_p$ there exists a solution for the equation $f_k(n)=p$, i.e. $N_k^*$ is well defined.

\item
We prove the statement for two separate cases: (1.) $s_k\rightarrow\infty$, and (2.) $s_k\rightarrow\bar{s}\in\mathbb{N}$. In the proof of both cases define the function $\upsilon$ as follows:	
\be
\upsilon(r):=\frac{\ln\ln(r)-\ln(4\pi)}{2\sqrt{2\ln(r)}}\ , \ \forall r>0
\ee 

{\bf \underline{Case $1$}:} \quad For any $\tau\in(0,1)$ define the following sequences:		
\be
N_k^{\pm}(\tau):=\frac{\sigma^2}{\Delta^2}\Big[\sqrt{(1\pm \tau)2\ln(k-s_k)}-\upsilon(k-s_k)+\sqrt{(1\pm \tau)2\ln(s_k)}-\upsilon(s_k)\Big]^2 .
\ee	
Recalling the normalizing sequences used in the proof of Statement $1$, basic limit arithmetics lead to the following limits:	
\be
L^*\big(N_k^{\pm}(\tau)\big) := \lim\limits_{k\rightarrow\infty}\Big\{\frac{\Delta\sqrt{N_k^{\pm}(\tau)}}{\sigma}-b_{k-s_k}-b_{s_k}\Big\}=\pm \infty . 
\ee	
As mentioned in the proof of Statement $1$, all the preconditions of Theorem \ref{thm:sum_max_a_k_positive} are satisfied. Therefore, Theorem \ref{thm:sum_max_a_k_positive} implies the next two limits for $f(n,k)$:
$$
\lim\limits_{k\rightarrow\infty}f\big(N^-_k(\tau),k\big)=0
$$
\be
\lim\limits_{k\rightarrow\infty}f\big(N^+_k(\tau),k\big)=1 .
\ee	
In addition, by definition, $\forall k>K_p$, $N^*_k$ satisfies $f(N^*_k,k)=p\in(0,1)$. Hence, $\exists K^1_\tau\geq K_p$ such that 
\be
f\big(N^-_k(\tau),k\big)< f\big(N^*_k,k\big)<f\big(N^+_k(\tau),k\big)\ \ , \ \ \forall k>K^1_\tau .
\ee	
$f(n,k)$ is nondecreasing in $n$. Therefore,
\be
N^-_k(\tau)\leq N^*_k\leq N^+_k(\tau)\ \ , \ \ \forall k>K^1_\tau . 
\label{eq:sample_size_sandwitch}		
\ee	
Recalling the exact expressions of $N^{\pm}_k(\tau)$, then 
\be
\lim\limits_{k\rightarrow\infty}\frac{N^{\pm}_k(\tau)}{\frac{2\sigma^2(1+\sqrt{C})^2}{\Delta^2}\ln(k-s_k)}= 1 \pm \tau . 
\label{eq:two_limits_scaled}
\ee	
Fix $\epsilon>0$. Since $k-s_k>s_k\rightarrow\infty$, the denominator of the two limits in eq. (\ref{eq:two_limits_scaled}) is positive for any large enough $k$. Therefore, eq. \eqref{eq:sample_size_sandwitch} implies that $\exists K_\tau^2\geq K_\tau^1$ such that 	
\be
1-\tau-\frac{\epsilon}{2}\leq\frac{N^*_k}{\frac{2\sigma^2(1+\sqrt{C})^2}{\Delta^2}\ln(k-s_k)}\leq 1+\tau+\frac{\epsilon}{2}\ \ , \ \ \forall k>K^2_{\tau} . 
\ee	
Thus, by rearranging this inequality and setting $\tau_\epsilon=\min\{\frac{\epsilon}{4},\frac{1}{2}\}\in(0,1)$, $\exists K_\epsilon:=K^2_{\tau_{\epsilon}}$ such that
\be
\Bigg|\frac{N^*_k}{\frac{2\Delta^2(1+\sqrt{C})^2}{\sigma^2}\ln(k-s_k)}-1\Bigg|\leq \min\big\{\frac{3\epsilon}{4},\frac{1}{2}+\frac{\epsilon}{2}\big\} <\epsilon\ \ , \ \ \forall k>K_{\epsilon}
\ee	
which ends the proof of this case. \\

{\bf \underline{Case $2$}:} \quad	The sequence of naturals $\{s_k;k\in\mathbb{N}\}$ satisfies $s_k \rightarrow\bar{s}\in\mathbb{N}$ hence $s_k=\bar{s}\in\mathbb{N}$ up to some finite prefix. Therefore, since this claim is about the asymptotic behavior of $\{N_k^*;k\in\mathbb{N}\}$ as $k\rightarrow\infty$, define the following sequences	
\be
N_k^{\pm}(\tau):=\frac{\sigma^2}{\Delta^2}\Big[\sqrt{(1 \pm \tau)2\ln(k-\bar{s})}-\upsilon(k-\bar{s})\Big]^2
\ee	
and observe that the following limits exist:	
\be
L^*\big(N_k^{\pm}(\tau)\big):=\lim\limits_{k\rightarrow\infty}\Big\{\frac{\Delta\sqrt{N_k^{\pm}(\tau)}}{\sigma}-b_{k-\bar{s}}\Big\}= \pm \infty . 
\ee
Therefore, Case $2$ is proven using exactly the same arguments used previously for Case $1$.
\end{enumerate}

\endproof

\subsection{Numerical Results}
\label{sec:one_stage_numeric}
\begin{figure}[!ht]
\begin{center}
\hspace*{-0.3in}
\includegraphics[height=2.85in, width=5.0in]{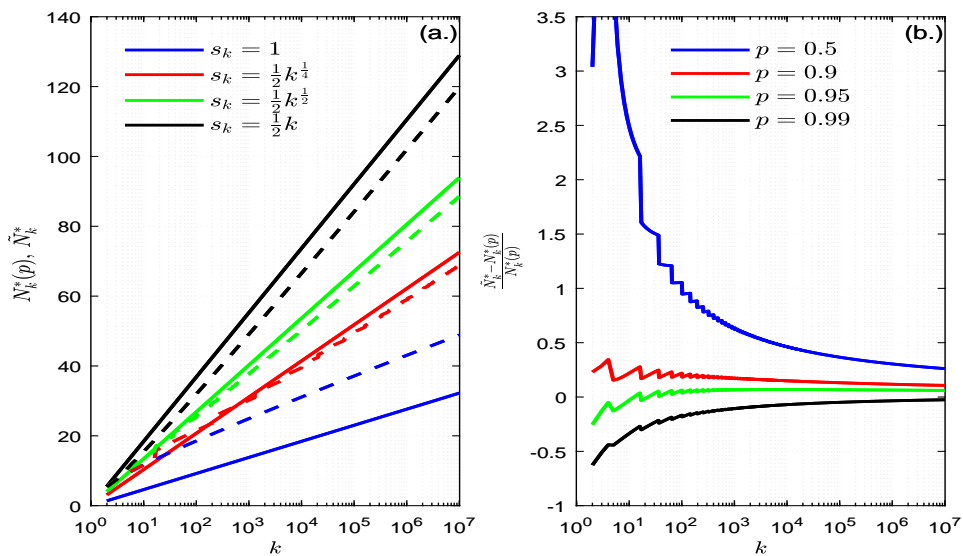} \hspace*{-0.0in} 
\caption{Asymptotic approximation of $N_k^*(p)$ for $\Delta=\sigma^2=1$ and for the {\it LFC} $\theta_{[i]} := \Delta \textbf{1}_{\{i>k-s_k\}}$. {\bf (a.)} The asymptotic approximation $\tilde{N}_k^*$ (dashed lines) and the exact sample size $N_k^*(0.95)$ (solid lines) vs. $k$ (on log-scale) for $s_k \propto k^{\alpha}$ for different values of $\alpha$. For all choices of $s_k$, the slopes of the $\tilde{N}_k^*$ lines, $2(1+\sqrt{\alpha})^2$, match the observed slope for true sample size $N_k^*(0.95)$, indicating that their ratio approaches $1$ as $k \to \infty$ {\bf (b.)} The relative error of the asymptotic approximation $\tilde{N}_k^*$ vs. the true sample size $N_k^*(p)$ for $s_k = \frac{1}{2} k^{\frac{1}{2}}$ for different values of $p$. While for small $k$ the probability of correct selection $p$ greatly affects the sample size, as $k \to \infty$ the sample size $N_k^*(p)$ becomes insensitive to $p$, and the asymptotic approximation $\tilde{N}_k^*$ gives the correct first order behavior for any fixed $p \in (0, 1)$.} 
\label{fig:one_stage_asymptotics}
\end{center}
\end{figure}

We next tested the accuracy of the asymptotic result in Theorem \ref{thm:generalized_sigmound_robbins} for finite $k$. To this end, we solved numerically eq. \eqref{eq:h_1_integral} 
to get the values of $N_k^*(p)$, and compared it to the asymptotic approximation $\tilde{N}_k^*$. 
Figure \ref{fig:one_stage_asymptotics} shows the approximation quality for different values of $p$, with specific values displayed in Table \ref{tab:one_stage_asymptotics}. The results confirm our analytic asymptotic predictions, yet the rate at which the numerical results approach the asymptotic limit depends on $p$ and the size $s_k$ of the selected set. 

\begin{table}[h!]
\begin{center}
\begin{tabular}{|c|c|c|c|c|} \hline 
$k$ \textbackslash\: $p$ & 0.5 & 0.9 & 0.95 & 0.99 \\ \hline
10 & 2.487 & 0.219 & -0.030 & -0.336 \\
100 & 1.055 & 0.219 & 0.064 & -0.167 \\ 
1000 & 0.637 & 0.175 & 0.069 & -0.105 \\ 
10000 & 0.461 & 0.149 & 0.069 & -0.071 \\ 
100000 & 0.367 & 0.132 & 0.067 & -0.049 \\ 
1000000 & 0.305 & 0.118 & 0.064 & -0.035 \\ 
10000000 & 0.261 & 0.107 & 0.061 & -0.025 \\ \hline 
\end{tabular} 
\caption{Relative error of asymptotic approximation $\frac{\tilde{N}_k^*-N_k^*(p)}{N_k^*(p)}$ for specific values in Figure \ref{fig:one_stage_asymptotics}.b}
\label{tab:one_stage_asymptotics}
\end{center}
\end{table}

\section{Two-Stage Procedures}
\label{sec:two_stage}
This Section revisits two well-known two-stage procedures which were suggested respectively in \cite{dudewicz1975allocation} and by \cite{rinott1978two}. Both procedures were developed for the problem of selecting the population with the largest mean from $k+1$ independent Gaussian populations with unknown and possibly different variances. Sections \ref{sec:stat_framework_two_stage}-\ref{sec:rinott} introduce the statistical settings and the two procedures. 
In Section \ref{sec:two_stage_asymptotics} we use Theorem \ref{thm:sum_max_a_k_zero} to analyze their asymptotic statistical efficiency as $k\rightarrow\infty$. Since both procedures draw a random number of samples, statistical efficiency is measured in terms of expected sample-size. 	
Our major conclusion states that as $k\rightarrow\infty$, the procedure in \cite{dudewicz1975allocation} is relatively more efficient by a factor of $2^{\frac{2}{N_0-1}}$, where $N_0$ is the sample size used by both procedures in the first stage. 

\subsection{Statistical Framework}
\label{sec:stat_framework_two_stage}
Let $X_{ij}\sim N(\theta_j,\sigma_j^2);i=1,\ldots,k+1, \: j\in\mathbb{N}$ be independent univariate Gaussian r.v's from the population $\Pi_i$ with unknown means $\theta_i$ and variances $\sigma_i^2$. Denote the ordered means by $\theta_{[1]}\leq\ldots\leq\theta_{[k+1]}$. The goal is to select the best population, namely the population whose mean is $\theta_{[k+1]}$. The settings here is similar to that of Section \ref{sec:robbins_siegmund}, except that the variances are unknown and may be different. We consider general selection procedures, namely multi-stage procedures which sequentially draw samples from the populations where the number of samples drawn from each population at any stage may depend on the sampling results of previous stages. 

In Section \ref{sec:robbins_siegmund} we analyzed a single-stage procedure. Considering the known indifference parameter $\Delta>0$ and the restricted parameter-space (see eq. (\ref{eq:indifference_parameter_space})) $\Theta(\Delta,k+1):=\big\{(\theta_1,\ldots,\theta_{k+1});\theta_{[k+1]}-\theta_{[k]} \geq \Delta\big\}$, 	
it was proven in \cite{dudewicz1971non} that no {\it single-stage} procedure controls the probability of correct selection above a prescribed value $p\in(0,1)$ for any parametrization $(\theta_1,\ldots,\theta_{k+1})\in\Theta(\Delta,k+1), (\sigma_1^2,\ldots,\sigma_{k+1}^2)\in\mathbb{R}^{k+1}_{++}$. 

Consequentially, \cite{dudewicz1975allocation,rinott1978two} provided two versions of {\it two-stage} procedures and have shown that they are guaranteed to control $\mathbb{P}(CS)$ above a prescribed value $p\in(0,1)$. We focus here on these two-stage procedures and describe them in details in the next subsections. 

\subsection{\DD's Procedure}
\DD \cite{dudewicz1975allocation} suggested a two-stage procedure $P_E$ which generalizes Stein' approach \cite{stein1945two}, described in Algorithm \ref{alg:dudewicz_dalal}. 

\begin{algorithm}[!h]
\caption{\DD's Two-Stage Procedure $P_E$}
\justify {\bf Input:} $\Delta$ - indifference parameter, $N_0 \geq 2$ - initial sample size, $p$ - desired $\mathbb{P}(CS)$ \\
{\bf Output:} $\hat{i}$ - selected population \\
{\bf Stage One:} 
\begin{enumerate} 
\item
Draw $N_0$ observations from each population $\Pi_i;i=1,\ldots,k+1$ 
\item
Compute the casual unbiased estimate $S^2_i$ for $\sigma^2_i$ from the initial sample taken for each $\Pi_i$	
\end{enumerate}
{\bf Stage Two:} 
\begin{enumerate} 
\item
For each $i=1,\ldots,k+1$ draw $N_i-N_0$ more samples from $\Pi_i$, where $N_i$ is given by
\be
N_i=\max\bigg\{N_0+1,\bigg \lceil \bigg(\frac{h_k^{(1)}}{\Delta}\bigg)^2 S_i^2\bigg\rceil \bigg\} .
\label{eq:sample_size_given_h}
\ee
Here $\lceil y \rceil$ denotes the smallest integer which is $\geq y$, and the constant $h_k^{(1)}$ is specified in eq. \eqref{eq:h1_integral_def}. 
\item Select numbers $\{a_{ij};j=1,\ldots,N_i\}_{i=1}^{k+1}$ such that $\forall i=1,\ldots,k+1$:

\begin{enumerate}
\item $S^2_i\sum_{j=1}^{N_i}a_{ij}^2=\bigg(\frac{\Delta}{h_k^{(1)}}\bigg)^2$ 
\item $\sum_{j=1}^{N_i}a_{ij}=1$ 
\item $a_{i1}=\ldots=a_{iN_0}$ 		
\end{enumerate}

\item Select the population $\hat{i}$ by the rule	
\be
\hat{i} := \underset{i=1,\ldots,k+1}{\arg\max} \big\{\sum_{j=1}^{N_i} a_{ij} X_{ij};\big\}
\ee
\end{enumerate}
\label{alg:dudewicz_dalal}
\end{algorithm}

Stage Two of Algorithm \ref{alg:dudewicz_dalal} requires calculating the numbers $\{a_{ij}\}$ and $h_k^{(1)}$. As mentioned in \cite{dudewicz1975allocation} a set of numbers $\{a_{ij}\}$ almost surely exists and it is easy to compute. 
The constant $h_k^{(1)}$ is chosen to guarantee a desired probability of correct selection $p$. This probability is bounded from below by the following integral	
\be
\mathbb{P}(CS|P_E)\geq\int_{-\infty}^{\infty}G_{\nu}^k(t+h_k^{(1)})g_{\nu}(t)dt
\ee	
where $G_{\nu}$ and $g_{\nu}$ are the c.d.f. and p.d.f. of Student's $t$-distribution with $\nu=N_0-1$ degrees of freedom (d.f's). Therefore, in order to ensure that the probability of correct selection remains above $p\in(0,1)$, $h_k^{(1)} := h_k^{(1)}(\nu)$ is determined as the solution of the following equation in $h$:
\be
f_1(h,k):=\int_{-\infty}^{\infty}G_{\nu}^k(t+h)g_{\nu}(t) dt = p .
\label{eq:h1_integral_def} 
\ee	
Although $h_k^{(1)}(\nu)$ depends on the initial sample size $N_0=\nu+1$, we usually omit $\nu$ and use the notation $h_k^{(1)}$ since $\nu$ is pre-defined and obvious from context. The next lemma ensures the validity of the asymptotic results in Subsection \ref{sec:two_stage_asymptotics}. 

\begin{lemma}
There exists $K>0$ such that $\forall k>K$ eq. \eqref{eq:h1_integral_def} has a unique positive solution. 
\label{lemma:unique_h}
\end{lemma}

\proof
$G_{\nu}$ is strictly increasing hence $f_1(h,k)$ is strictly increasing in $h$, $\forall k\in\mathbb{N}$. Since $G_{\nu}$ is bounded and satisfies $G(-\infty)=0, G_{\nu}(\infty)=1$, by the bounded convergence theorem $f_1(-\infty,k)=0$, $f_1(\infty,k)=1$. Similarly, the continuity of $f_1(\cdot,k)$ on $\mathbb{R}$ can be also justified by the bounded convergence theorem and hence by the intermediate value theorem $\forall k\in\mathbb{N}$ there exists unique $h_k^{(1)}\in\mathbb{R}$ such that $f_1(h,k)=p$. In addition, $G_{\nu}^k(t)\xrightarrow{k\rightarrow\infty}0, \: \forall t\in\mathbb{R}$ and hence by the bounded convergence theorem, $f_1(0,k)\xrightarrow{k\rightarrow\infty}0<p$. Therefore, since $f_1(h,k)$ is strictly increasing in $h$, deduce that $\exists K\in\mathbb{N}$ such that $f_1(h,k)<p, \: \forall h\leq 0, \forall k>K$. Finally, because $\forall k>K, \: \exists h_k^{(1)}$ which satisfies $f_1(h_k^{(1)},k)=p$, the constant $h_k^{(1)}$ must be positive $\forall k>K$. 
\endproof

\subsection{Rinott's Procedure}
\label{sec:rinott}
Since the procedure $P_E$ allows some means to be negatively weighted, Rinott \cite{rinott1978two} stated that as pointed in \cite{stein1945two}, a similar procedure based on ordinary means may be more appealing. Rinott introduced such a procedure $P_R$ which guarantees a probability of correct selection above $p$ and shares the same steps of $P_E$ except two differences: First, in Step $3$ set $a_{ij}=\frac{1}{N_i}$, $\forall i,j$; second, in Step $2$ replace $h_k^{(1)}$ by another well-defined sequence $h_k^*\geq h_k^{(1)}$ which is determined as a solution of a certain integral equation specified by Rinott. Practically, Rinott suggested to use another sequence $ h_k^{(2)} := h_k^{(2)}(\nu) \geq h^*_k$ which is defined for each $k\in\mathbb{N}$ as the solution of the following simpler equation	
\be
f_2(h,k):=\Big[\int_{-\infty}^{\infty}G_{\nu}(t+h)g_{\nu}(t)dt\Big]^k=p . 
\label{eq:h2_integral_def}
\ee	
The same arguments used in Lemma \ref{lemma:unique_h} for $\{h_k^{(1)}\}_{k\in\mathbb{N}}$ show that $\{h_k^{(2)}\}_{k\in\mathbb{N}}$ is also well-defined and positive up to a finite prefix. Consequently, since our analysis performs an asymptotic comparison of the procedures $P_E$ and $P_R$, w.l.o.g. we make the simplifying assumption that $0<h_k^{(1)}\leq h_k^{(2)}, \: \forall k\in\mathbb{N}$. 

\subsection{Asymptotic Efficiency}
\label{sec:two_stage_asymptotics}	
Since both of the procedures depicted previously a draw random number of samples, it is convenient to determine their asymptotic efficiency by the expected sample-size required in order to satisfy the $\mathbb{P}(CS)$ criterion as $k\rightarrow\infty$. To see how this expected sample size relates to $h_k^{(j)},\: j=1,2$, observe that both procedures draw an infinite number of samples as $k\rightarrow\infty$. Therefore, regardless the value of $N_0$, for any large enough $k$, the procedures $P_E$ and $P_R$ are associated respectively with expected sample-sizes of $\big(h_k^{(1)}\big)^2 \sum_{i=1}^{k+1}\frac{\sigma_i^2}{\Delta^2}$ and $\big(h_k^{(2)}\big)^2 \sum_{i=1}^{k+1}\frac{\sigma_i^2}{\Delta^2}$, up to rounding errors. Consequently, in order to analyze the asymptotic relative efficiency in terms of expected sample-size, it is enough to determine the asymptotic behavior of the ratio $\Big(\frac{h_k^{(2)}}{h_k^{(1)}}\Big)^2$ as $k\rightarrow\infty$. Theorems \eqref{thm:dalal_asymptotic} and \eqref{thm:rinott_procedure} make the first order approximations $h_k^{(1)} \sim \tilde{h}_k^{(1)}; h_k^{(2)}\sim \tilde{h}_k^{(2)}$ as $k\rightarrow\infty$ with $\tilde{h}_k^{(2)} := 2^{\frac{1}{\nu}}\tilde{h}_k^{(1)} := 2^{\frac{1}{\nu}}\gamma_{\nu} k^{\frac{1}{\nu}}q_p$, where $q_p$ is the $p$'th quantile of $\nu$-Fr\'echet distribution and $\gamma_{\nu}$ is some function of $\nu$ to be specified later. This result implies that for any initial sample size $h_k^{(2)}-h_k^{(1)}\rightarrow\infty$. Therefore, regardless the exact value of $p$, the numerical insight made in the last paragraph of Subsection 4.1 of \cite{rinott1978two} which states that for $p \geq 0.75$, the difference between $h_k^{(1)}$ and $h_k^{(2)}$ should be small, is incorrect for large enough values of $k$ unless the sample size $N_0=\nu+1$ is also increased.

\begin{theorem}
\label{thm:dalal_asymptotic}
Let $q_p$ be the $p$'th quantile of $\nu$-Fr\'echet distribution and let $\gamma_{\nu}$ be defined as follows		
\be
\gamma_{\nu} := \Bigg[\frac{\Gamma(\frac{\nu+1}{2})}{\nu^{1-\frac{\nu}{2}} \sqrt{\pi}\Gamma(\frac{\nu}{2})}\Bigg]^{\frac{1}{\nu}} .
\ee		
Then $h_k^{(1)}\sim \tilde{h}_k^{(1)} := \gamma_\nu k^{\frac{1}{\nu}}q_p$.		
\end{theorem}

\proof	
Set some $\tau\in(0,p\wedge1-p)$ and define the following sequences:	
\be
h^{\pm}_k(\tau):=\gamma_\nu k^{\frac{1}{\nu}}q_{p \pm \tau}
\ee		
where $\gamma_{\nu}$ has already been defined in the statement of the theorem and $q_{p\pm\tau}$ are the $p\pm\tau$'th quantiles of the $\nu$-Fr\'echet distribution. 

The Fr\'echet distribution is nonnegative and continuous and hence its quantiles are simply defined by the inverse of the $\nu$-Fr\'echet c.d.f ${e^{-x}}^{-\nu}, \: \forall x>0$, for any $\tau>0$ such that $0<p-\tau<p+\tau<1$. In addition, let $X_1,\ldots,X_{k+1}\overset{i.i.d}{\sim}t_\nu$ be a sequence of $k+1$ i.i.d Student's $t$ r.v's with $\nu$ d.f's. By the convolution formula for difference of independent r.v's, $f_1(h,k)$ can be expressed as follows:	
\be		
f_1(h,k)=\int_{-\infty}^{\infty}G_{\nu}^k(t+h)g_{\nu}(t)dt=\mathbb{P}\Big(\max_{j=1,\ldots,k}X_j-X_{k+1}\leq h\Big) .
\ee		
Recall a known result (see e.g. Proposition $2.5$ in \cite{grigelionis2013student}, with the constant $\gamma_{\nu}$ corrected here) which states that the extreme value distribution of a sequence of i.i.d Student's $t$ random variables with $\nu$ d.f's is $\nu$-Fr\'echet distribution with the normalizing constants $a_k:=\gamma_{\nu}^{-1} k^{-\frac{1}{\nu}}$ and $b_k\equiv0$. Denote the following limits		
\be
L^{**}\big(h_k^{\pm}(\tau)\big)=\lim\limits_{k\rightarrow\infty}\gamma_\nu^{-1} k^{-\frac{1}{\nu}} \gamma_\nu k^{\frac{1}{\nu}}q_{p \pm \tau}=q_{p \pm \tau}. 
\ee			
The c.d.f of $\nu$-Fr\'echet distribution is continuous on $\mathbb{R}$ and in particular on $\{q_{p-\tau},q_{p+\tau}\}$. Therefore, Theorem \ref{thm:sum_max_a_k_zero} can be used to obtain that	
$$		
\lim\limits_{k\rightarrow\infty}f_1\big(h_k^+(\tau),k\big)=p+\tau>p
$$
\be
\lim\limits_{k\rightarrow\infty}f_1\big(h_k^-(\tau),k\big)=p-\tau<p
\ee		
By definition, $f_1(h_k^{(1)},k)=p, \: \forall k\in\mathbb{N}$ and hence by simple limit rules $\exists K_\tau\in\mathbb{N}$ such that 	
\be		
f_1\big(h_k^-(\tau),k\big)<f_1(h_k^{(1)},k)<f_1\big(h_k^+(\tau),k\big)\ \ , \ \ \forall k>K_\tau .
\ee		
Since for any $k\in\mathbb{N}$, $f_1(h,k)$ is strictly increasing in $h$, then 	
\be		
h_k^-(\tau)<h_k^{(1)}< h_k^+(\tau)\ , \ \forall k>K_\tau . 
\label{eq:h_1_sandwitch}
\ee		
The Fre\'chet distribution is nonnegative and continuous and hence $q_p>0$, i.e. $\forall k\in\mathbb{N}$, $\gamma_\nu k^\frac{1}{\nu}q_p>0$. Thus, dividing eq. (\ref{eq:h_1_sandwitch}) by $\gamma_{\nu} k^{\frac{1}{\nu}}q_p$ gives	
\be		
\frac{q_{p-\tau}}{q_p}\leq \frac{h_k^{(1)}}{\gamma_{\nu} k^\frac{1}{\nu}q_p}\leq \frac{q_{p+\tau}}{q_p}\ , \ \forall k>K_\tau . 
\ee		
Since the $\nu$-Fr\'echet distribution is continuous, the inverse function theorem states that $q_{p+\tau}$ and $q_{p-\tau}$ are continuous functions of $\tau$ on $(0,p \wedge 1-p)$. Therefore, by taking $\tau\rightarrow0^+$ both boundaries approach to 1, i.e. $\forall \epsilon>0$, $\exists \tau_\epsilon\in(0,p\wedge1-p)$ which satisfies	
\be		
\max\Big\{|\frac{q_{p-\tau_\epsilon}}{q_p}-1|,|\frac{q_{p+\tau_\epsilon}}{q_p}-1|\Big\}<\epsilon
\ee	
and respectively	
\be		
-\epsilon<\frac{h_k^{(1)}}{\gamma_\nu k^\frac{1}{\nu}q_p}-1<\epsilon\ , \ \forall k>K_{\tau_\epsilon}
\ee		
which, by definition, is an equivalent writing of the needed result $h_k^{(1)}\sim \tilde{h}_k^{(1)} := \gamma_\nu k^\frac{1}{\nu}q_p$. 

\endproof

\begin{theorem}
\label{thm:rinott_procedure}
With the same notations of Theorem \ref{thm:dalal_asymptotic}, $h_k^{(2)}\sim \tilde{h}_k^{(2)} := 2^{\frac{1}{\nu}}\gamma_\nu k^{\frac{1}{\nu}}q_p$.
\end{theorem}

\proof
Let $T_i^{(j)}\overset{i.i.d}{\sim}t_\nu ; i=1,\ldots,k ; j=1,2 $ be Student's $t$ r.v's . Due to the symmetry of Student's $t$-distribution around zero, the convolution formula for difference of independent r.v's implies that $f_2(h,k)$ can be expressed as follows:	
\begin{align}
f_2(h,k) &:= \Big[\int_{-\infty}^{\infty}G_\nu(t+h)g_\nu(t)dt\Big]^k \nonumber \\
&= \mathbb{P}\big(\max_{i=1,\ldots,k}\{T_i^{(1)}-T_i^{(2)}\}\leq h\big) \nonumber \\		
&= \mathbb{P}\big(\max_{i=1,\ldots,k}\{T_i^{(1)}+T_i^{(2)}\}\leq h\big) . 
\label{eq:rinott_integral_max}
\end{align}
Let $\tilde{g}_\nu$ be the density associated with the distribution of a sum of two i.i.d Student's $t_\nu$ r.v's. This density is given in eq. (2.1) of \cite{ghosh1975distribution}
\be
\tilde{g}_{\nu}(t) = \frac{\Gamma(\frac{\nu+1}{2}) \Gamma(\nu+\frac{1}{2})}{2^{\nu}\nu^{\frac{1}{2}} \Gamma^2(\frac{\nu}{2}) \Gamma(\frac{\nu}{2}+1)} \Big(\frac{4\nu}{4\nu+t^2}\Big)^{1+\nu} \hygeom\Big(\frac{1}{2}, \nu+\frac{1}{2}; \frac{\nu}{2}+1; \frac{t^2}{4\nu+t^2}\Big)
\ee
where $\hygeom(a, b; c; z)$ is the hypergeometric function with parameters $a,b,c$ evaluated at $z$ with $|z|<1$. We next use Euler's transformation for the hypergeometric function, 
\be
\hygeom(a, b; c; z) = (1-z)^{c-a-b}\hygeom(c-a, c-b; c; z)
\ee
to get:
\be
\tilde{g}_{\nu}(t) = \frac{\Gamma(\frac{\nu+1}{2}) \Gamma(\nu+\frac{1}{2})}{2^{\nu}\nu^{\frac{1}{2}} \Gamma^2(\frac{\nu}{2}) \Gamma(\frac{\nu}{2}+1)} \Big(\frac{4\nu}{4\nu+t^2}\Big)^{1+\frac{\nu}{2}} \hygeom\Big(\frac{\nu+1}{2}, \frac{1-\nu}{2}; \frac{\nu}{2}+1; \frac{t^2}{4\nu+t^2}\Big). 
\label{eq:student_sum_density_hypergeometric}
\ee
We have $\frac{t^2}{4\nu+t^2} \underset{t \to \infty}{\longrightarrow} 1$ therefore, 
\begin{align}
\hygeom\Big(\frac{\nu+1}{2}, \frac{1-\nu}{2}; \frac{\nu}{2}+1; \frac{t^2}{4\nu+t^2}\Big) &\underset{t \to \infty}{\longrightarrow} \hygeom\Big(\frac{\nu+1}{2}, \frac{1-\nu}{2}; \frac{\nu}{2}+1; 1\Big) \nonumber \\
&= \frac{\Gamma(\frac{\nu}{2}+1)\Gamma(\frac{\nu}{2})}{\Gamma(\frac{1}{2})\Gamma(\nu+\frac{1}{2})}
\label{eq:Gauss_hypergeometric_1}
\end{align}
where the value of the hypergeometric function evaluated at $1$ is an analytical continuation which is provided by Gauss' Theorem. Plugging eq. \eqref{eq:Gauss_hypergeometric_1} into eq. \eqref{eq:student_sum_density_hypergeometric} and taking $t \to \infty$ we get: 
\be
\tilde{g}_{\nu}(t) \sim \frac{2 \Gamma(\frac{\nu+1}{2})}{\nu^{\frac{-\nu}{2}}\sqrt{\pi} \Gamma(\frac{\nu}{2})} t^{-(1+\nu)}\sim 2 g_{\nu}(t) .
\ee
Since the asymptotic values of the densities $\tilde{g}_{\nu}(t), {g}_{\nu}(t)$ for large $t$ are the same up to a multiplicative factor of $2$, we can follow Propositions $2.3$ of \cite{grigelionis2013student} to get the asymptotic cumulative distribution function of $\tilde{g}_{\nu}(t)$ for $t \to \infty$: 
\be
1 - \tilde{G}_{\nu}(t) \sim \frac{2 \Gamma(\frac{\nu+1}{2})}{\nu^{1-\frac{\nu}{2}}\sqrt{\pi} \Gamma(\frac{\nu}{2})} t^{-\nu} = 2 \gamma_{\nu}^{\nu} t^{-\nu} \sim 2 \big(1-{G}_{\nu}(t)\big)
\ee
and follow Propositions $2.5$ of \cite{grigelionis2013student} to get 
that the extreme value distribution for $\tilde{g}_{\nu}(t)$
is the $\nu$-Fr\'echet distribution with the normalizing constants $a_k:=2^{-\frac{1}{\nu}}\gamma_{\nu}^{-1} k^{-\frac{1}{\nu}}$ and $b_k\equiv0$.

Next, set some $\tau\in(0,p\wedge1-p)$ and define the following sequences:
\be
h^{\pm}_k(\tau)=a_k^{-1}q_{p \pm \tau} . 
\ee	
Finally, the arguments used to prove Theorem \ref{thm:dalal_asymptotic} hold for this case too and hence imply the needed result.
\endproof

\begin{corollary}
\label{corr_proc_difference}
Using the same notations of Theorem \ref{thm:dalal_asymptotic}, $h_k^*-h_k^{(1)}=\mathcal{O}(\gamma_\nu k^{\frac{1}{\nu}}q_p)$ as $k\rightarrow\infty$.
\end{corollary}

\proof
In \cite{rinott1978two}, it was shown that $\forall k\in\mathbb{N}$, $h_k^{(1)}\leq h_k^*\leq h_k^{(2)}$, therefore $\forall k\in\mathbb{N}$:
\be
0\leq\frac{h_k^*-h_k^{(1)}}{\gamma_\nu k^{\frac{1}{\nu}}q_p}\leq \frac{h_k^{(2)}-h_k^{(1)}}{\gamma_\nu k^{\frac{1}{\nu}}q_p}\xrightarrow{k\rightarrow\infty} 2^{\frac{1}{\nu}}-1 .
\ee	
This limit implies that 	
\be
0\leq\limsup_{k\rightarrow\infty}\frac{h_k^*-h_k^{(1)}}{\gamma_\nu k^{\frac{1}{\nu}}q_p}\leq 2^{\frac{1}{\nu}}-1<\infty
\ee	
and the corollary follows from the definition of the big $\mathcal{O}$ notation.
\endproof

Theorems \ref{thm:dalal_asymptotic},\ref{thm:rinott_procedure} show the dependency of the asymptotics of $h_k^{(1)}(\nu)$, $h_k^{(2)}(\nu)$ on the initial sample size $N_0=\nu+1$. In particular, the asymptotic relative efficiency of the two procedures satisfies 
\be
\lim_{\nu \to \infty} \lim_{k \to \infty} \Big(\frac{h_k^{(2)}(\nu)}{h_k^{(1)}(\nu)}\Big)^2 = \lim_{\nu \to \infty} 2^{\frac{2}{\nu}} = 1 . 
\label{eq:limit_nu_k}
\ee
The next theorem reveals that surprisingly, when we fix first $\nu=\infty$ and then let $k\rightarrow\infty$, the limit is given by: 
\be
\lim_{k \to \infty} \Big(\frac{h_k^{(2)}(\infty)}{h_k^{(1)}(\infty)}\Big)^2 = 2 
\label{eq:limit_k_nu}
\ee
i.e. $\nu=\infty$ is a discontinuity point of the asymptotic relative efficiency as a function of $\nu$.

\begin{theorem}
For $\nu=\infty$ we have $h_k^{(2)}(\infty) \sim\sqrt{2}h_k^{(1)}(\infty) \sim 2\sqrt{\ln(k)}$
\label{thm:ratio_Gaussians}
\end{theorem}

\proof 
\begin{enumerate}
\item To derive the asymptotics of $h_k^{(2)}(\infty)$, recall that a standard Student's $t$-distribution with $\nu=\infty$ d.f's is a standard Gaussian distribution. Thus, let $X_1,X_2,\ldots\stackrel{i.i.d.}{\sim} N(0,2)$ and observe that 

\be
f_2(h,k):=\mathbb{P} \Big(\max_{j=1,\ldots,k} X_j\leq h\Big)=\mathbb{P}\Big(\max_{j=1,\ldots,k}Z_j\leq\frac{h}{\sqrt{2}}\Big)
\ee

where $Z_1,Z_2,\ldots\stackrel{i.i.d.}{\sim}N(0,1)$. As mentioned in Section \ref{sec:robbins_siegmund}, the normalizing constants of the standard Gaussian distribution are $a_k=\sqrt{2\ln(k)}$ and $b_k\sim\sqrt{2\ln(k)}$. Thus, for any $\tau\in(0,p\wedge 1-p)$ define $h_k(\tau)^+:=\sqrt{2}(b_k-\frac{g_{p\pm\tau}}{a_k})$ where $g_{p\pm\tau}$ is the $p \pm\tau$'s quantile of the standard Gumbel distribution. Since the extreme value distribution of the standard Gaussian distribution is a standard Gumbel distribution, this implies that
\be
\lim_{k\rightarrow\infty} f_2 \big(h_k^{\pm}(\tau),k\big)=p\pm\tau .
\ee
Therefore, by the same technique which was used in the proof of Theorem \ref{thm:dalal_asymptotic}, deduce that $h_k^{(2)}(\infty) \sim \sqrt{2}(b_k-\frac{q_p}{a_k})\sim\sqrt{2}b_k \sim 2 \sqrt{\ln(k)}$.

\item To derive the asymptotics of $h_k^{(1)}$, set an arbitrary $\tau\in{(0,p\wedge 1-p)}$ and define $h_k^\pm(\tau)=b_k-z_{p\pm\tau}$ where $z_p$ is the $p$'th quantile of the standard Gaussian distribution. Let $Z_1,Z_2,\ldots \stackrel{i.i.d}{\sim}N(0,1)$. For
\be
f_1(h,k):=\mathbb{P}\Big(\max_{j=1,\ldots,k} Z_j+Z_{k+1}\leq h\Big)
\ee
Theorem \ref{thm:sum_max_a_k_positive} implies that
\be
\lim_{k\rightarrow\infty}f\big(h_k^\pm(\tau),k\big)=p\pm\tau . 
\ee
Thus, the same technique used to prove Theorem \ref{thm:dalal_asymptotic} shows that $h_k^{(1)}(\infty) \sim b_k \sim \sqrt{2\ln(k)}$. 
\end{enumerate}

\endproof

\begin{remark}
Theorems \ref{thm:dalal_asymptotic}-\ref{thm:ratio_Gaussians} state that the relative asymptotic efficiency of the procedures is invariant to the value of $p$. 
\end{remark}

\begin{remark}
For both procedures, the guaranteed lower bounds on the probability of correct selection may not be tight and hence an empirical comparison of sample sizes giving the same probability of correct selection in practice may give different conclusions and should be studied separately. This issue was studied using simulations in \cite{branke2007selecting,wang2013reducing}. Lately, \cite{Frazier2014} gave new theoretical insights regarding this phenomenon. 	
\end{remark}

\subsection{Numerical Results}
\label{sec:two_stage_numeric}
\begin{figure}[!ht]
\begin{center}
\hspace*{-0.2in} \includegraphics[height=5.8in, width=5.0in]{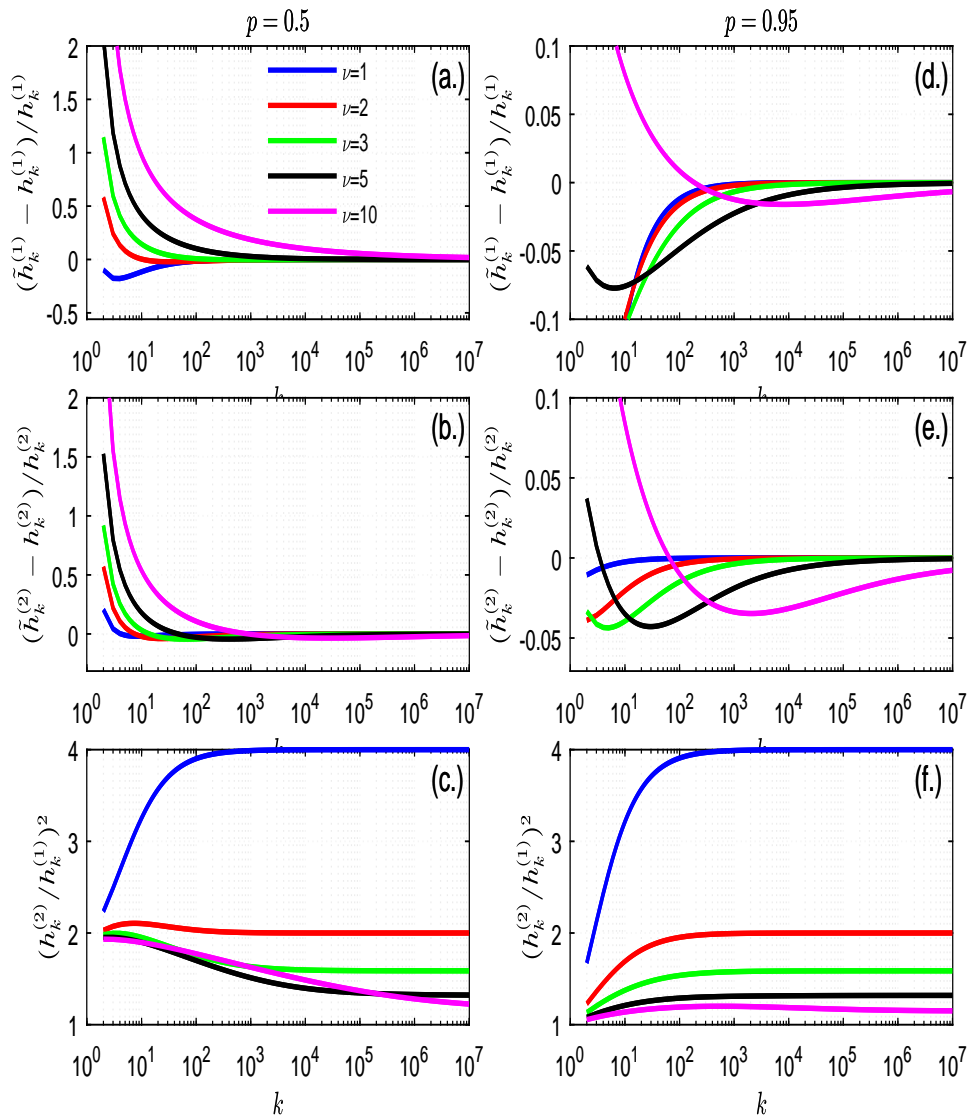} \vspace*{-0.2in} 
\caption{Comparison of relative efficiency of two procedures for $p=0.5,0.95$ for different values of $\nu$.
	{\bf (a.)} The relative error of the asymptotic approximation $\tilde{h}_k^{(1)}$ vs. the exact (numeric) integral $h_k^{(1)}$ for \DD's procedure. For $k \sim 10^2-10^4$ the relative error of the asymptotic approximation is only a few percents; the approximation accuracy decreases when $\nu$ is increased. {\bf (b.)} The relative error of the asymptotic approximation $\tilde{h}_k^{(2)}$ vs. the exact (numeric) integral $h_k^{(2)}$ for Rinott's procedure. The qualitative behavior of the approximation is similar to that of $h_k^{(1)}$, with larger relative error as $\nu$ is increased. {\bf (c.)} The relative efficiency of the two procedures approaches the asymptotic value of $2^{\frac{2}{\nu}}$ as $k \to \infty$. For larger values of $\nu$, larger $k$ values are needed to get an accurate approximation. {\bf (d.)-(f.)} The same as (a.)-(c.) but for $p=0.95$, showing improved accuracy for the asymptotic approximation is for larger $p$. The approximation accuracy for $\tilde{h}_k^{(1)}$ here is not monotonic with $\nu$.} 
\label{fig:two_stage_relative}
\end{center}
\end{figure}

We solved numerically eq. \eqref{eq:h1_integral_def} and \eqref{eq:h2_integral_def}
to get the values of $h_k^{(1)}$ and $h_k^{(2)}$, respectively, and compared them with the asymptotic results in Theorems \ref{thm:dalal_asymptotic},\ref{thm:rinott_procedure} for finite $k$.
Figure \ref{fig:two_stage_relative} shows the relative efficiency of the two procedures, with specific values displayed in Table \ref{tab:two_stage_relative}. The results confirm our asymptotic predictions. The rate at which the numerical results approach the asymptotic limit varies with $\nu$ and $p$. 

\begin{table}[h!]
\begin{center}
\begin{tabular}{|c|c|c|c|c|} \hline
$k$ & $\frac{\tilde{h}_k^{(1)}-h_k^{(1)}}{h_k^{(1)}}, p=0.5$ & $\frac{\tilde{h}_k^{(2)}-h_k^{(2)}}{h_k^{(2)}}, p=0.5$ & $\frac{\tilde{h}_k^{(1)}-h_k^{(1)}}{h_k^{(1)}}, p=0.95$ & $\frac{\tilde{h}_k^{(2)}-h_k^{(2)}}{h_k^{(2)}}, p=0.95$ \\ \hline
10 & 0.975 & 0.537 & 0.079 & 0.084 \\ 
100 & 0.375 & 0.107 & 0.009 & -0.010 \\
1000 & 0.186 & -0.003 & -0.012 & -0.034 \\ 
10000 & 0.101 & -0.033 & -0.016 & -0.032 \\ 
100000 & 0.056 & -0.033 & -0.014 & -0.022 \\ 
1000000 & 0.032 & -0.023 & -0.010 & -0.013 \\ 
10000000 & 0.018 & -0.014 & -0.007 & -0.008 \\ \hline 
\end{tabular}

\caption{Relative error of asmymptotic approximation for two procedures $(\tilde{h}_k^{(i)}-h_k^{(i)}) / h_k^{(i)}$ for specific values in Figure \ref{fig:two_stage_relative}}
\label{tab:two_stage_relative}
\end{center}
\end{table}

\newpage
\clearpage

\subsection{Choosing the parameter $\nu$}
The statistical efficiency of the two procedures in \cite{dudewicz1975allocation,rinott1978two} depend on the choice of the parameter $\nu$. In this section we derive for each procedure an asymptotic approximation for the optimal $\nu$ minimizing the expected sample size as $k\rightarrow\infty$. Define the expected sample size for the two procedures when choosing the parameter $\nu$, 
\be
\mu_k^{(j)} := \mu_k^{(j)}(\nu) = \mathbb{E} [{\bf N_k^{(j)}}] = \sum_{i=1}^{k+1} \mathbb{E} [N_i], \: j=1,2
\label{eq:expected_sample_size}
\ee
where $N_{i,k}=\max\bigg(N_0+1, \Big(\frac{\sigma_i h_k^{(j)}}{\Delta}\Big)^2 \bigg)$ are given in eq. \eqref{eq:sample_size_given_h}, and ${\bf N_k^{(j)}}= \sum_{i=1}^{k+1} N_{i,k}$ is the actual (random) sample size. Since $\mu_k^{(j)}(\nu) \to \infty$ as $\nu \to \infty$, a minimizer $\nu_k^{(j*)}$ must exist. Thus, we may define the optimal parameter choice and the optimal sample size attained for the two procedures 
\be
\nu_k^{(j*)} := \arg\min_{\nu \in \mathbb{N}} \mu_k^{(j)}(\nu) \: ; \: \mu_k^{(j*)} := \min_{\nu \in \mathbb{N}} \mu_k^{(j)}(\nu) = \mu_k^{(j)}(\nu_k^{(j*)}) . 
\label{eq:optimal_expected_sample_size}
\ee

Finding the optimal parameter $\nu_k^{(j*)}$ leads to both conceptual and technical difficulties. First, the optimum depends on the unknown variances $\sigma_i^2$. Second, even if the variances $\sigma_i^2$ were known, the maximization operation and the non-explicit form of $h_k^{(j)}$ makes the optimization problem computationally challenging. 

To overcome these difficulties, we propose a parameter choice for $\nu$ based on two simplifications: 
(i) We ignore the maximization with $N_0+1$ in the definition of $N_{i,k}$ and optimize only the second term as we take $k \to \infty$, and (ii) we replace $h_k^{(j)}$ by its asymptotic approximation $\tilde{h}_k^{(j)}$. With these two simplifications, we define the approximate expected sample size
\be
{\tilde{\mu}}_k^{(j)}(\nu) := \big(\tilde{h}_k^{(j)}(\nu)\big)^2 \frac{\sum_{i=1}^k \sigma_i^2}{\Delta^2}, \: j=1,2 \: .
\label{eq:approx_expected_sample_size}
\ee
and the approximate optimal parameter
\be
{\tilde{\nu}_k}^{(j*)} := \arg\min_{\nu \in \mathbb{N}} {\tilde{\mu}}_k^{(j)}(\nu) \: ; \: {\tilde{\mu}}_k^{(j*)} := \min_{\nu \in \mathbb{N}} \tilde{\mu}_k^{(j)}(\nu) = \tilde{\mu}_k^{(j)}({\tilde{\nu}_k}^{(j*)} ) . 
\label{eq:approx_optimal_expected_sample_size}
\ee
${\tilde{\nu}_k}^{(j*)}, {\tilde{\mu}}_k^{(j*)}$ do not depend on the unknown variances $\sigma_i^2$ and can be found by minimizing: 
\be
\tilde{h}_k^{(1)}(\nu) = \gamma_v k^{\frac{1}{\nu}} q_p = \Bigg[\frac{\Gamma(\frac{\nu+1}{2}) k }{-\nu^{1-\frac{\nu}{2}} \sqrt{\pi}\Gamma(\frac{\nu}{2}) \ln(p)}\Bigg]^{\frac{1}{\nu}} .
\label{eq:h_nu}
\ee

\begin{theorem}
For $k$ large enough, eq. \eqref{eq:h_nu} has a unique solution ${\tilde{\nu}_k}^{(j*)}$. Moreover, as $k \to \infty$: ${\tilde{\nu}_k}^{(j*)} \sim 2 \ln(k)$, $\tilde{h}_k^{(1*)} \sim \sqrt{2 e \ln(k)}$ and ${\tilde{\mu}}_k^{(j*)} \sim 2 e \ln(k) \frac{\sum_{i=1}^k \sigma_i^2}{\Delta^2}$ . 
\label{thm:optimal_asymptotic}
\end{theorem}

\proof
Differentiating the logarithm of eq. \eqref{eq:h_nu}, 
\be
\ln \big(\tilde{h}_k^{(1)}(\nu)\big) = \frac{1}{\nu} \Bigg[ \ln\big(\frac{-k}{\sqrt{\pi} \ln(p)}\big) + \ln \big(\Gamma(\frac{\nu+1}{2})\big) - \ln\big(\Gamma(\frac{\nu}{2})\big) + (\frac{\nu}{2}-1) \ln(\nu) \Bigg]
\label{eq:h_nu_log}
\ee
gives the first order condition:
\be
0 = \frac{d \ln \tilde{h}_k^{(1)}(\nu)}{d \nu} = \frac{1}{2 \nu^2} \mathcal{H}_k(\nu)
\label{eq:first_order_nu}
\ee
where 
\be
\mathcal{H}_k(\nu) := -2 -2 \ln\big(\frac{-k}{\sqrt{\pi} \ln(p)}\big) + \nu + 2 \ln(\nu) + 2 \ln \Big(\frac{\Gamma(\frac{\nu}{2})}{\Gamma(\frac{\nu+1}{2})}\Big) + \nu \Big( \Psi\big(\frac{\nu+1}{2}\big) - \Psi\big(\frac{\nu}{2}\big)\Big) 
\label{eq:nu_first_order}
\ee
and $\Psi$ is the digamma function. Since $\nu>0$ we have $sign\big(\frac{d \ln \tilde{h}_k^{(1)}(\nu)}{d \nu}\big) = sign\big(\mathcal{H}_k(\nu)\big)$ and the first order condition is satisfied if and only if $\mathcal{H}_k(\nu) = 0$. The derivative of $\mathcal{H}_k$ is 
\be
\frac{d \mathcal{H}_k(\nu)}{d \nu} = 1 + \frac{2}{\nu} + \frac{\nu}{2} \Big(\Psi'(\frac{\nu+1}{2}) - \Psi'(\frac{\nu}{2})\Big) . 
\label{eq:h_derivative}
\ee
By Lemma $1$ in \cite{alzer2001mean}, $\Psi'$ is strictly monotonically decreasing in $\mathbb{R}_{+}$. Therefore, using the recurrence relation for polygamma functions $\Psi'(z+1)=\Psi'(z)-\frac{1}{z^2}$ we get the bound 
\be
\Psi'(\frac{\nu+1}{2}) > \Psi'(\frac{\nu}{2}+1) = \Psi'(\frac{\nu}{2}) - \frac{4}{\nu^2}.
\label{eq:polygamma_recurrence}
\ee
Plugging eq. \eqref{eq:polygamma_recurrence} into eq. \eqref{eq:h_derivative} gives $\frac{d \mathcal{H}_k(\nu)}{d \nu} > 1, \quad \forall \nu>0$, hence $\mathcal{H}_k$ is monotonically increasing. We use the following bounds, 
\begin{align}
-2 \ln(\nu) &< 2 \ln \Big(\frac{\Gamma(\frac{\nu}{2})}{\Gamma(\frac{\nu+1}{2})}\Big) < 0 \nonumber \\
0 &< \nu \Big( \Psi\big(\frac{\nu+1}{2}\big) - \Psi\big(\frac{\nu}{2}\big)\Big) < 1 
\end{align}
to bound $\mathcal{H}_k(\nu)$
\be
-2 -2 \ln\big(\frac{-k}{\sqrt{\pi} \ln(p)}\big) + \nu < \mathcal{H}_k(\nu) < -1 -2 \ln\big(\frac{-k}{\sqrt{\pi} \ln(p)}\big) + \nu + 2 \ln(\nu) . 
\label{eq:h_derivative_sandwitch}
\ee
For $k$ large enough the right bound in eq. \eqref{eq:h_derivative_sandwitch} shows that $\mathcal{H}_k(1) < 0$. 
For any fixed $k$, eq. \eqref{eq:h_derivative_sandwitch} gives $\mathcal{H}_k(\nu) \sim \nu \to \infty$ as $\nu \to \infty$. Since $\mathcal{H}_k(\nu)$ is monotonically increasing for $\nu>0$, the first order condition in eq. \eqref{eq:nu_first_order} has a unique solution which is the global minimum of $\ln \big(\tilde{h}_k^{(1)}(\nu)\big)$ in $\nu \in (1,\infty)$. 

We can solve eq. \eqref{eq:nu_first_order} numerically to get the optimal $\nu$ for any given $k$ and $p$. To get the asymptotic solution as $k \to \infty$ we set $\mathcal{H}_k(\nu)=0$ in eq. \eqref{eq:h_derivative_sandwitch}, 
\begin{align}
-2 -2 \ln(\sqrt{\pi}\ln(p)) + \nu &< 2 \ln(k) < -1 -2 \ln(\sqrt{\pi}\ln(p)) + \nu + 2 \ln(\nu) \nonumber \\
&\implies {\tilde{\nu}_k}^{(1*)} \sim 2 \ln(k) . 
\label{eq:nu_first_order_asymptotics}
\end{align}
Plugging the asymptotic solution ${\tilde{\nu}_k}^{(1*)} \sim 2 \ln(k)$ into the asymptotic expression for $h_k^{(1)}$ yields
\begin{align}
\tilde{h}_k^{(1)}\big(2 \ln (k)\big) &= \gamma_{2\ln(k)} k^{\frac{1}{2\ln(k)}} q_p \nonumber \\
&= \Bigg[\frac{\Gamma\big(\ln(k)+\frac{1}{2}\big) k }{-\big(2\ln(k)\big)^{1-\ln(k)} \sqrt{\pi}\Gamma\big(\ln(k)\big) \ln(p)}\Bigg]^{\frac{1}{2\ln(k)}} \nonumber \\
&\sim \sqrt{2e\ln(k)} . 
\end{align}
Thus, for the choice ${\tilde{\nu}_k}^{(j*)}= 2 \ln(k)$ the approximate expected asymptotic sample size for \DD's procedure is $\tilde{\mu}_k^{(1*)} \sim 2 e \ln(k) \frac{\sum_{i=1}^k \sigma_i^2}{\Delta^2}$. 

\endproof

\begin{remark}
It is instructive to compare Theorem \ref{thm:optimal_asymptotic} in the case of equal variances $\sigma_i \equiv \sigma$ to \RS's one-stage procedure applied when the variance is known. \RS's procedure \cite{robbins1967iterated} requires $\sim 2 \ln(k) \frac{\sigma^2}{\Delta^2}$ samples from each population in order to ensure correct selection with a prescribed probability $p$, i.e. the overall sample size summing over all populations is $\sim 2 k\ln(k) \frac{\sigma^2}{\Delta^2}$. Hence the approximate asymptotic sample size for the case of unknown variance is within a multiplicative factor of $e$ of the sample size for the case of known variance. 
\end{remark}

For Rinott's procedure \cite{rinott1978two}, the asymptotic behavior of $\tilde{h}_k^{(2)}$ can be similarly derived, yielding 
\be
\tilde{h}_k^{(2)}\big(2 \ln (k)\big)\sim 2^{\frac{1}{2\ln(k)}} \tilde{h}_k^{(1)}\big(2 \ln (k)\big) \sim \tilde{h}_k^{(1)}\big(2 \ln (k)\big) \sim \sqrt{2e\ln(k)}. 
\ee
Thus, while for every fixed $\nu$ \DD's procedure is asymptotically more efficient, as $k\rightarrow\infty$, the approximations of the optimal $\nu$'s for both procedures are equivalent up to a first order error term. The reason is that as $\nu = 2\ln(k)$ increases, the asymptotic ratio $2^{\frac{1}{2\ln(k)}}$ goes to $1$. Although to the first order the two sample sizes are identical, taking a multiplicative factor $2^{\frac{1}{2\ln(k)}}$ into account for Rinott's procedure may yield more accurate results. 

We next study the asymptotic behavior of $h_k^{(2)}$ for fixed $k$ as $\nu \to \infty$. Lemma \ref{lemma:monotonic_h} shows a useful monotonicity property of Student's $t$ r.v's, which is 
used to show the monotonicity of $h_k^{(2)}$ in $\nu$. 
\begin{lemma}
Let $T_i^{(j)} \sim t_{\nu_i} \: ; i,j=1,2$ be four independent Student's $t$ r.v's. 
with ${\nu}_1\leq{\nu}_2$. Then $\forall h>0$
\be
\mathbb{P}(T_1^{(1)} + T_1^{(2)} \leq h) \leq \mathbb{P}(T_2^{(1)} + T_2^{(2)} \leq h) .
\ee
\label{lemma:monotonic_h}
\end{lemma}

\proof
Let $\nu_1 < \nu_2$. Using the symmetry of the Student's $t$ densities $g_{\nu_j}; j=1,2$ around zero, we get 
\begin{align}
\mathbb{P}(T_2^{(1)} + T_2^{(2)} \leq h) &= \int_{-\infty}^{\infty} G_{\nu_2}(h-t) g_{\nu_2}(t) dt \nonumber \\
&= \int_{-\infty}^{\infty} G_{\nu_2}(t) g_{\nu_2}(t-h) dt \nonumber \\ 
&= \int_{0}^{\infty} G_{\nu_2}(t) g_{\nu_2}(t-h) dt + \int_{0}^{\infty} G_{\nu_2}(-t) g_{\nu_2}(t+h) dt \nonumber \\ 
&= \int_{0}^{\infty} G_{\nu_2}(t) [g_{\nu_2}(t-h)-g_{\nu_2}(t+h)] dt + 1 - G_{\nu_2}(h) \nonumber \\
&\geq \int_{0}^{\infty} G_{\nu_1}(t) [g_{\nu_2}(t-h)-g_{\nu_2}(t+h)] dt + 1 - G_{\nu_2}(h) \nonumber \\
&= \int_{-\infty}^{\infty} G_{\nu_1}(h-t) g_{\nu_2}(t) dt \nonumber \\
&= \int_{-\infty}^{\infty} G_{\nu_2}(t) g_{\nu_1}(t-h) dt \nonumber \\
&\geq \int_{-\infty}^{\infty} G_{\nu_1}(h-t) g_{\nu_1}(t) dt \nonumber \\
&= \mathbb{P}(T_1^{(1)} + T_1^{(2)} \leq h) 
\label{eq:rinott_monotonic_h}
\end{align}
where $G_{\nu_1}(t) < G_{\nu_2}(t), \:\: \forall t > 0$ and $g_{\nu_2}(t-h)-g_{\nu_2}(t+h) \geq 0, \:\: \forall t,h>0$ together imply the first inequality appearing in the fifth line of eq. \eqref{eq:rinott_monotonic_h} above. To obtain the second inequality, repeat lines $2$-$6$ 
of eq. \eqref{eq:rinott_monotonic_h} with $g_{\nu_2}$ replaced by $g_{\nu_1}$.
\endproof

\begin{corollary}
$\exists K>0$ such that $\forall k>K$ , $h_k^{(2)}(\nu)$ is monotonically non-increasing in $\nu$. 
\label{cor:monotonic_h}
\end{corollary}
\proof
By Lemma \ref{lemma:unique_h}, $\exists K>0$ such that $\forall k>K$, $h_k^{(2)}(\nu_1), h_k^{(2)}(\nu_2)>0$. Using eq. \eqref{eq:h2_integral_def} for $h=h_k^{(2)}(\nu_2)>0$ and representing the probabilities in Lemma \ref{lemma:monotonic_h} as convolutions,
\be
\int_{-\infty}^{\infty} G_{\nu_1}(h_k^{(2)}(\nu_2)-t) g_{\nu_1}(t) dt \leq
\int_{-\infty}^{\infty} G_{\nu_2}(h_k^{(2)}(\nu_2)-t) g_{\nu_2}(t) dt = p^{\frac{1}{k}} 
\ee
and because $G_{\nu_1}$ is monotonically increasing we get $h_k^{(2)}(\nu_1) \geq h_k^{(2)}(\nu_2)$. 
\endproof

\begin{figure}[!ht]
\begin{center}
\hspace*{-0.52in} \includegraphics[height=2.5in, width=5.5in]{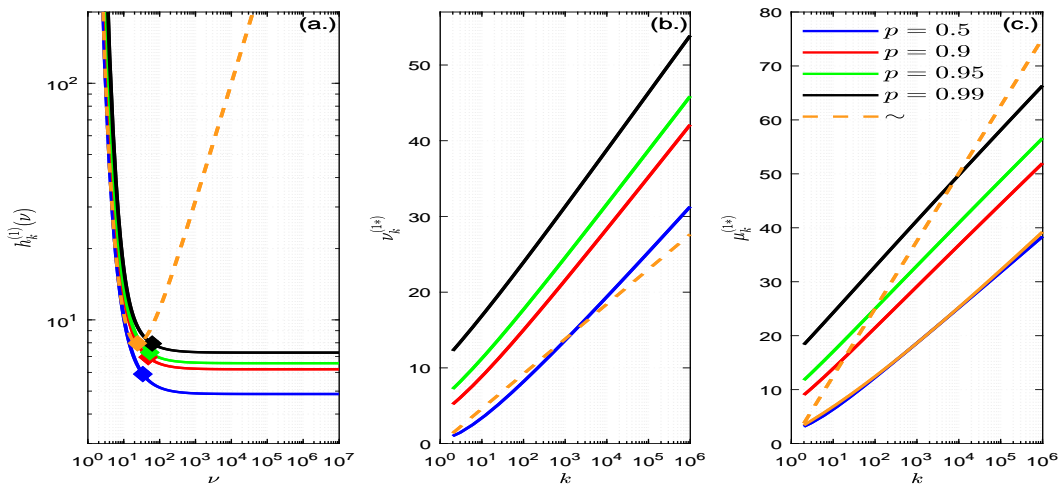} \vspace*{-0.1in} 
\caption{Optimal choice of parameters for \DD's procedure for different values of $p$ and $k$, for $\sigma_i^2 \equiv \Delta^2 = 1$. 
	{\bf (a.)} The exact $h_k^{(1)}(\nu)$ for different values of $p$ (solid lines of different colors) vs. the asymptotic $\tilde{h}_k^{(1)}(\nu)$ (dashed orange line), for $k=10^5$, as a function of $\nu$. The orange diamond represent the minimum $(\tilde{\nu}_k^{(1*)}, \tilde{h}_k^{(1*)})$ as found by the first order condition in eq. \eqref{eq:first_order_nu}. 
	The other colored diamonds are the optimal values $(\nu_k^{(1*)}, h_k^{(1*)})$ minimizing $h_k^{(1)}(\nu)^2$ for different values of $p$. In all cases the computed $h_k^{(1)}(\nu)^2$ was monotonically decreasing in $\nu$, and $\nu_k^{(1*)}$ was the solution of $\nu+2 = h_k^{(1)}(\nu)^2$. 
	{\bf (b.)} The optimal parameter $\nu_k^{(1*)}$ for different values of $p$, and the approximate optimal value $\tilde{\nu}_k^{(1*)}$, as a function of $k$, shown on a log-scale. The slope for the approximation is lower than the slope of the exact lines, indicating that for large $k$, $\tilde{\nu}_k^{(1*)}$ underestimate $\nu_k^{(1*)}$ for this case. 
	{\bf (c.)} The resulting optimal expected sample size $\mu_k^{(1*)}$ as function of $k$. For the approximation the approximate sample size $\tilde{\mu}_k^{(1*)}$ (dashed orange) and the exact expected sample size for $p=\frac{1}{2}$, evaluated at the approximate optimum, $\mu_k^{(1)}(\tilde{\nu}_k^{(1*)})$ (solid orange) are shown. The slope for the approximation $\tilde{\mu}_k^{(1*)}$ is higher than the slope of the exact expected sample sizes. However, the exact expected sample size at our approximate solution $\mu_k^{(1)}(\tilde{\nu}_k^{(1*)})$ matches for this case the optimal exact expected sample size for $p=\frac{1}{2}$. } 
\label{fig:optimal_nu}
\end{center}
\end{figure}

\begin{remark}
The order of limits of $k,\nu \to \infty$ matters. Eq. (\ref{eq:limit_nu_k}) shows that for fixed finite $\nu$, 
$\tilde{h}_k^{(j)}(\nu) \sim h_k^{(j)}(\nu)$ as $k \to \infty$. Now, for fixed $k$,
$\tilde{h}_k^{(j)}(\nu)$ and $h_k^{(j)}(\nu)$ behave differently as $\nu \to \infty$. The approximation $\tilde{h}_k^{(j)}(\nu)$ has a unique minimum at $(1, \infty)$ as we have shown above, and deriving the asymptotics of eq. \eqref{eq:h_nu} shows 
\be
\tilde{h}_k^{(1)}(\nu) \sim \frac{1}{\nu^{(\frac{1}{\nu}-1)/2}}\sim\sqrt{\nu}\xrightarrow{\nu\rightarrow\infty}\infty .
\label{eq:tilde_h_asymptotics}
\ee
In contrast, for $k$ fixed and large enough, Lemma \ref{lemma:monotonic_h} implies that the exact $h_k^{(2)}(\nu)$ for Rinott's procedure is monotonically decreasing in $(1, \infty)$ such that $\underset{\nu \to \infty}{\lim} h_k^{(2)}(\nu) = \frac{\Phi^{-1}(p^{\frac{1}{k}})}{\sqrt{2}}$. In addition, since $h_k^{(1)}(\nu) \leq h_k^{(2)}(\nu)$ by Proposition $3$ in \cite{rinott1978two}, $h_k^{(1)}(\nu)$ is bounded from above by a monotonically decreasing sequence, and $\underset{\nu\rightarrow\infty}{\limsup}\: h_k^{(1)}(\nu) \leq \frac{\Phi^{-1}(p^{\frac{1}{k}})}{\sqrt{2}}$. Thus, the asymptotic convergence $\underset{k \to \infty}{\lim} \frac{h_k^{(j)}}{\tilde{h}_k^{(j)}} = 1$ shown in Theorems \ref{thm:dalal_asymptotic},\ref{thm:rinott_procedure} is therefore not uniform in $\nu$. In particular, as shown in Figure \ref{fig:optimal_nu}, the asymptotic result in Theorem \ref{thm:optimal_asymptotic} does not necessarily imply $h_k^{(1)}\big(2 \ln (k)\big) \sim \sqrt{2e\ln(k)}$ and $\mu_k^{(1)}(2 \ln(k)) \sim 2 e \ln(k) \frac{\sum_{i=1}^k \sigma_i^2}{\Delta^2}$. 
\end{remark}

Finally, recall that our simplification defined $\tilde{\mu}_k^{(j)}$ as an approximate expected sample size, while ignoring the maximization in the definition of $N_{i,k}$. In Theorem \ref{thm:asymptotic_optimality_nu} we define an approximate sample size which does take the maximization into account and show optimality with respect to this definition. Since for bounded sequences $\{\nu_k\}$ with $\underset{k \to \infty}{\limsup} \:\nu_k = \nu \in \mathbb{N}$ Theorem \ref{thm:dalal_asymptotic} shows $\mu_k^{(1)}(\nu_k) = \Omega(k^{\frac{1}{\nu}}) \gg \ln(k)$, we consider only sequences $\{\nu_k\}$ such that $\nu_k\rightarrow \infty$.
For any such sequence we give conditions ensuring that, almost surely for each population the asymptotic approximate sample size as $k\rightarrow\infty$: (i) converges to its expectation, and (ii) cannot be improved compared to $\{\tilde{\nu}_k^{(1*)} \}$.

\newpage
\clearpage

\begin{theorem}
Consider a sequence $\{\nu_k\}$ such that $\nu_k\rightarrow\infty$ and for each $k\in\mathbb{N}$, let 

\begin{enumerate}
\item $\tilde{N}_{i,k}:=\max\big\{\nu_k+2,\big[\tilde{h}_k^{(1)}(\nu_k)\big]^2\frac{S_i^2}{\Delta^2}\big\}\ , \ \ \forall i=1,\ldots,k+1$. 

\item $\tilde{N}_{i,k*}:=\max\big\{\tilde{\nu}_k^{(i*)}+2,\big[\tilde{h}_k^{(1)}(\tilde{\nu}_k^{(i*)})\big]^2\frac{S_i^2}{\Delta^2}\big\} \ , \ \ \forall i=1,\ldots,k+1$.

\end{enumerate}

If $2e \sigma_i^2 > \Delta^2, \: \forall i\in\mathbb{N} $, then 
\begin{enumerate}

\item $\tilde{N}_{i,k*}\sim\big[\tilde{h}_k^{(1)}(\tilde{\nu}_k^{(i*)})\big]^2\frac{\sigma_i^2}{\Delta^2} \ \ , \ \ \forall i\in\mathbb{N} \ \ , \ \ \mathbb{P}-a.s.$

\item $\underset{k\rightarrow\infty}{\liminf}\:\frac{\tilde{N}ֶֶֶ_{i,k}}{\tilde{N}_{i,k*}} \geq 1\ , \ \ \forall i\in\mathbb{N}\ , \ \ \mathbb{P}-a.s.$

\end{enumerate}
\label{thm:asymptotic_optimality_nu}
\end{theorem}

\proof
Since $\nu_k\rightarrow\infty$ and the $S_i^2$'s are unbiased estimators of $\sigma_i^2$'s computed using independent samples from independent populations, we may use the strong law of large numbers throughout the proof. 

\begin{enumerate}

\item 
The following approximation stems from the strong law of large numbers and eq. \eqref{eq:tilde_h_asymptotics}:
\be
\big[\tilde{h}_k^{(1)}(\nu_k)\big]^2\frac{S_i^2}{\Delta^2}\sim \big[\tilde{h}_k^{(1)}(\nu_k)\big]^2 \frac{\sigma_i^2}{\Delta^2}\ , \ \ \forall i\ \ , \ \ \mathbb{P}-a.s.
\ee
Let $i\in\mathbb{N}$ and consider the sequence $\{\tilde{\nu}_k^{(j*)}\}$. Theorem \ref{thm:optimal_asymptotic} states that $\tilde{h}_k^{(1)}(\tilde{\nu}_k^{(i*)})\sim\sqrt{2e\ln(k)}$ hence the fact that $2e\frac{\sigma_i^2}{\Delta^2}>1$ implies that almost surely, up to a finite prefix 
\be
\tilde{\nu}_k^{(i*)}+2 \leq \big[\tilde{h}_k^{(1)}(\tilde{\nu}_k^{(i*)})\big]^2\frac{S_i^2}{\Delta^2} .
\ee

Thus, almost surely, up to a finite prefix $\tilde{N}_{i,k*}$ is given by its second argument, i.e.
\be
\tilde{N}_{i,k*}\sim\big[\tilde{h}_k^{(1)}(\tilde{\nu}_k^{(i*)})\big]^2\frac{\sigma_i^2}{\Delta^2} \ , \ \mathbb{P}-a.s.
\ee
Therefore, since the intersection of a countable number of events of probability one is a an event of probability one, the claim follows.

\item Let $i\in\mathbb{N}$ and consider the following two cases: \\
(i) If the inequality $\nu_k\leq\big[\tilde{h}_k^{(1)}(\nu_k)\big]^2\frac{\sigma_i^2}{\Delta^2}$ holds for $k$ large enough, then the strong law of large numbers implies that $\tilde{N}_{i,k}\sim\big[\tilde{h}_k^{(1)}(\nu_k)\big]^2\frac{\sigma_i^2}{\Delta^2}$ and hence:
\be
\liminf_{k\rightarrow\infty}\frac{\tilde{N}ֶֶֶ_{i,k}}{\tilde{N}_{i,k*}}=\liminf_{k\rightarrow\infty}\frac{[\tilde{h}_k^{(1)}(\nu_k)\big]^2\frac{\sigma_i^2}{\Delta^2}}{[\tilde{h}_k^{(1)}(\nu_k^*)\big]^2\frac{\sigma_i^2}{\Delta^2}}\geq 1 . 
\ee
\end{enumerate}

\endproof

\newpage
\clearpage
\section{Discussion}
\label{sec:discussion}
In this work we proved limit theorems for linear combinations of partial maxima, and demonstrated their utility by using them to derive the asymptotic behaviors of several selection procedures under different statistical frameworks. The specific contributions to the R\&S theory in Sections \ref{sec:robbins_siegmund} and \ref{sec:two_stage} shed new light on existing popular procedures, and offer natural avenues for future research. In particular, Section \ref{sec:robbins_siegmund} studies a new asymptotic regime where the number of populations to be selected is determined as a function of the total number of populations. Studying the behavior of other R\&S procedures in this regime can lead to similar generalization of other known results. In Section \ref{sec:two_stage} we have shown that the guarantees of the procedure of \DD \cite{dudewicz1975allocation} are asymptotically superior to that of Rinott \cite{rinott1978two}. While several authors proposed new procedures based on Rinott's procedure \cite{ahmed2002simulation,boesel2003using,nelson2001simple}, it would be interesting to develop R\&S procedures based on \DD's procedure. If such new procedures are found, would they be better than the current Rinott's-based procedures? More questions related to the comparison between the two procedures in \cite{dudewicz1975allocation,rinott1978two} are

\begin{enumerate}
\item Can one apply the techniques we used for studying $h_k^{(1)},h_k^{(2)}$ to derive similar asymptotic results for $h_k^*$? 

\item Can one prove rigorously that our conjecture that $N_0 \sim 2 \ln(k)$ samples for the first stage of both procedures is optimal, hence the relative efficiency of the two procedures under optimal choice of $\nu$ is one? 
can similar results be proven for the asymptotic expected sample size of the two procedures under optimal choice of the $N_0$ parameter? 

\item What is the relative efficiency of the two procedures when considering the {\it actual} probability for correct selection instead of its bounds? 

\end{enumerate} 

Beyond the procedures discussed in this work, it would be interesting to apply our approach more generally to study the asymptotic attributes of other, more modern, R\&S procedures such as the ones proposed in \cite{Frazier2014,kim2001fully,nelson2001simple}.

Other asymptotic regimes for R\&S procedures can also be studied using tools from extreme value theory - for example, in \cite{mukhopadhyay1979some} it was shown that Rinott's procedure is asymptotically inefficient in the sense of \cite{chow1965asymptotic}, i.e. in the asymptotic regime where $\Delta^* \downarrow 0$. As the number of items $k$ is increased, it is of interest to study the case where the indifference parameter $\Delta^*$ is decreased, for example $\Delta^*(k) \propto k^{-1}$. 
This case arises naturally when the populations have parameters $\theta_i$ within a prespecified range, or drawn from a certain prior distribution in a Bayesian setting, as was studied for example in \cite{ein2006thousands,zuk2007ranking}. As the number of selected items $s_k$ is increased, it is also of interest to relax the requirement for correct selection, and allow approximate correct selection, for example requiring correct selection of $(1-\delta) s_k$ items for some predefined $\delta>0$. 

Taking a broader view, this work points to an interesting relation between extreme value theory (a nice introduction is provided by the books \cite{deoliveira1979asymptotic,leadbetter2012extremes}) and the asymptotic behavior of R\&S procedures. Therefore, other results from extreme value theory can be naturally applied to R\&S procedures - for example, it would be interesting yet challenging to develop and apply limit theorems for maxima of dependent random variables, in order to study R\&S procedures for dependent populations. 

Finally, the limit theorems proved in Section \ref{sec:linear_maxima} may be applied to other fields beyond that of R\&S procedures. In a well known application of extreme value theory, it is used to calculate the statistical significance of a local sequence alignment in computational biology \cite{dembo1994limit}. In this application, deriving the distribution of the best (maximal) sequence alignment under the null is required in order to establish whether two aligned sub-sequences are significantly similar, in an hypothesis testing framework. Sometimes a single sequence alignment does not provide sufficient statistical evidence against the null, and pooling information from several local sequence alignments in the same region is required - hence the need to calculate the distribution of the sum of several maxima under the null, which can hopefully be achieved using our theorems. 
We hope that the current work will stimulate search for further applications and generalizations of our theorems.

\subsection*{Acknowledgements}
The authors express their deep gratitude to Pavel Chigansky for fruitful conversations and to Yosef Rinott for sharing his insights about his previous work.

\newpage
\clearpage
\bibliographystyle{apalike}
\bibliography{ranking_selection}
\end{document}